%% file: main.tex
\documentclass[draft,onecolumn]{IEEEtran}
\usepackage[utf8]{inputenc}
\usepackage{amsmath,amsfonts}
\usepackage{amssymb}
\usepackage{bbold}
\usepackage{graphicx}
\usepackage{physics}
\usepackage{caption}
\usepackage{subcaption}
\usepackage{enumerate}
\usepackage{tikz}

\usepackage{cite} 
\usepackage{amsmath} 
\usepackage{amssymb}  
\usepackage{mathrsfs}
\usepackage{color}
\usepackage{nicefrac}
\usepackage[linesnumbered, ruled, vlined]{algorithm2e}
\usepackage{mathtools}
\mathtoolsset{showonlyrefs}
\let\IEEEproof\proof                         
\let\proof\relax

\usepackage{amsthm}

\usepackage{wrapfig,graphicx}
\usepackage{hyperref}
\graphicspath{{.}}

\theoremstyle{plain}

\newtheorem{problem}{Problem}
\usepackage{etoolbox}
\usepackage{float}
\usepackage{tabularx}

\def\th@definition{
	\thm@headfont{\itshape} 
	\thm@notefont{} 
}
\makeatother

\usepackage{color}

\theoremstyle{definition}
\newtheorem{definition}{Definition}

\makeatletter
\newenvironment{subtheorem}[1]{%
  \def\subtheoremcounter{#1}%
  \refstepcounter{#1}%
  \protected@edef\theparentnumber{\csname the#1\endcsname}%
  \setcounter{parentnumber}{\value{#1}}%
  \setcounter{#1}{0}%
  \expandafter\def\csname the#1\endcsname{\theparentnumber.\Alph{#1}}%
  \ignorespaces
}{%
  \setcounter{\subtheoremcounter}{\value{parentnumber}}%
  \ignorespacesafterend
}
\makeatother
\newcounter{parentnumber}

\newtheorem{lemma}{Lemma}
\newtheorem{example}{Example}
\newtheorem{theorem}{Theorem}
\newtheorem{remark}{Remark}

\newtheorem{corollary}{Corollary}

\newcommand{\R}{\mathbb{R}}
\newcommand{\N}{\mathbb{N}}

\newcommand{\adj}{\textnormal{Adj}}

\author{Calvin Hawkins, Bo Chen, Kasra Yazdani, Matthew Hale$^{\ast}$
\thanks{$^{\ast}$Department of  Mechanical and Aerospace Engineering at the University of Florida, Gainesville, FL USA. Emails: \texttt{\{calvin.hawkins,bo.chen,kasra.yazdani,matthewhale\}} \texttt{@ufl.edu}. This work was supported in part by NSF under CAREER Grant~{\#}1943275, by AFOSR under Grant~{\#}FA9550-19-1-0169,
and by ONR under Grant~{\#}N00014-21-1-2502. 
}
}

\title{Node and Edge Differential  Privacy  for Graph Laplacian Spectra: Mechanisms and Scaling Laws}

\date{November 2022}

\begin{document}

\maketitle
\begin{abstract}
    This paper develops a framework for privatizing the spectrum of the Laplacian of an undirected graph using differential privacy. 
    We consider two privacy formulations. The first obfuscates the presence of edges in the graph and the second obfuscates the presence of nodes. 
    We compare these two privacy formulations and show that the privacy formulation that considers edges is better suited to most engineering applications. We use the bounded Laplace mechanism to provide $(\epsilon,\delta)$-differential privacy to the eigenvalues of a graph Laplacian, and we pay special attention to the algebraic connectivity, 
    which is the Laplacian's the second smallest eigenvalue. Analytical bounds are presented on the accuracy of the mechanisms and on certain graph properties computed
    with private spectra. A suite of numerical examples confirms the accuracy of private spectra in practice. 
\end{abstract}
\input{introduction}
\input{background}
\input{privacy_mechanisms}

\input{scaling}
\input{other_properties}
\input{guidelines_examples}

\section{Conclusions}
This paper presented two differential privacy mechanisms for edge and node privacy of the
spectra of graph Laplacians of unweighted, undirected graphs. 
Bounded noise was used to provide private values that are still
accurate, and the private values of Laplacian spectrum
were shown to give accurate estimates of the diameter and mean distance
of a graph, the trace of the Laplacian, the Kemeny constant, and Cheeger's inequality. Future work includes the development of
new privacy mechanisms for other algebraic graph properties.

\bibliographystyle{IEEEtran}{}
\bibliography{sources}

\input{appendix}
\end{document}

%% file: introduction.tex
\section{Introduction}

Graphs are used to model a wide range
of interconnected systems, including multi-agent control
systems~\cite{Ren2005}, social networks~\cite{Scott1988}, and others~\cite{SHIRLEY2005287}.
Various properties of these graphs have been used to
analyze controllers and dynamical processes over them,
such as reaching a consensus~\cite{zheng_consensus_2011}, 
the spread of a virus~\cite{Mieghem2009}, robustness to 
connection failures~\cite{freitas2020evaluating}, and others.
Graphs in these applications may contain sensitive information,
e.g., one's close friendships in the case of a social network,
and it is essential that these analyses do not inadvertently
leak any such information.

Unfortunately, it is well-established that even graph-level analyses
may inadvertently reveal sensitive information about individuals in them,
such as the absence or presence of individual nodes in a graph~\cite{kasiviswanathan13}
and the absence or presence of specific edges between them~\cite{karwa14}.
Similar privacy threats have received attention in the
data science community, where graphs represent datasets
and the goal is to enable data analysis while safeguarding
the data of individuals in those datasets.

Differential privacy is one well-studied tool for doing so.
Differential privacy is a statistical notion of privacy that has several
desirable properties: (i) it is robust to side information, in that learning
additional information about data-producing entities does not weaken privacy
by much~\cite{kasiviswanathan14}, and (ii) it is immune to post-processing, in that arbitrary
post-hoc computations on private data do not weaken privacy~\cite{dwork_algorithmic_2013}. 
There exist numerous differential privacy implementations for
graph properties, including counts
of subgraphs~\cite{karwa14}, 
degree distributions~\cite{day16}, and 
other frequent patterns in graphs~\cite{shen13}. 
These privacy mechanisms generally follow
the pattern of computing the quantity of interest, adding
carefully calibrated noise to it, and releasing its
noisy form. Although simple, this approach strongly
protects data with a suite of guarantees provided
by differential privacy~\cite{dwork_algorithmic_2013}. 

The need for privacy for the aforementioned graph properties comes from the inferences that one can draw about a graph from these quantities, as detailed in~\cite{Ding2018, Hay2009, Task2012}. Decades of research in algebraic graph theory have quantified connections between the Laplacian spectrum and a myriad of other graph properties; see~\cite{abreu07} for a summary. Accordingly, the Laplacian spectrum, especially the algebraic connectivity $\lambda_2,$ implicates the same ability to draw inferences as other graph properties and hence gives rise to the same types of privacy concerns.

We therefore protect the values the graph Laplacian spectrum using two notions of privacy:
edge and node differential privacy \cite{karwa2011private}. Edge privacy
obfuscates the absence and/or presence of a pre-specified number
of edges, while node privacy obfuscates the absence or presence of a single node. In this paper we show that the differences in guarantees of these two notions of privacy result in drastic differences in the accuracy of the private values of the Laplacian spectrum. Specifically, in Section~\ref{sec:scaling} we show that the variance of noise required to obfuscate the presence of one node in a graph of size $n$ scales with $n^2$, which rapidly grows large. For this reason, Sections~\ref{sec:bo_section} and~\ref{sec:examples} focus on edge privacy and obfuscating the connections in a network.
We note that while differential privacy has been applied to protect various quantities in multi-agent systems~\cite{hawkins2020differentially,gohari2020privacy,gohari2021differential,cortes2016differential}, privacy for properties
of a multi-agent network itself has received less attention, and that is what we focus on. 

In this paper we pay special attention to the algebraic connectivity. A graph's algebraic connectivity (also called its Fiedler value~\cite{fiedler_algebraic_1973})
is equal to the second-smallest eigenvalue of its Laplacian.
This value plays a central role in the study of multi-agent
systems because it sets the convergence rates
of consensus algorithms~\cite{olfati04}, which appear directly or in modified
form in formation control~\cite{ren07}, 
connectivity control~\cite{gennaro06}, 
and many distributed optimization 
algorithms~\cite{nedic18}. 

Our implementation uses the recent bounded Laplace
mechanism~\cite{holohan2018bounded}, which  ensures that private scalars lie in a specified interval. The algebraic connectivity of a graph
is bounded below by zero and above by the number of nodes
in a graph, and we confine private outputs
to this interval by applying the mechanism in~\cite{holohan2018bounded} to the privatization of Laplacian spectra.

\emph{Contributions:} We provide closed-form values for the sensitivity
and other constants needed to define edge and node differential privacy mechanisms
for the Laplacian spectrum, and this is the first contribution
of this paper. The second contribution is showing the detrimental scaling of node privacy and the benefits of edge privacy.
Our third contribution is the use of the private values
of algebraic connectivity to analytically bound other graph properties, namely
the diameter of graphs and the mean distance between their nodes. 
Our fourth contribution is providing guidelines on using these mechanisms  
by providing a series of examples to demonstrate how to use the mechanisms and the accuracy of information they provide. 

We note that~\cite{wang13} has developed a different
approach to privacy for the eigendecomposition of
a graph's adjacency matrix. Given our motivation
by multi-agent systems, we focus
on a graph's Laplacian, which commonly appears in multi-agent controllers, and we derive
simpler forms for the distribution of noise required,
as well as a privacy mechanism that does not
require any post-processing. 

A preliminary version of this paper appeared in \cite{chen2021edge}. This paper extends the edge privacy mechanism for $\lambda_2$ to the rest of the Laplacian spectrum, develops the node privacy mechanism for $\lambda_2,$ compares the scaling of the edge and node privacy mechanisms, and provides further applications and uses of the private Laplacian spectrum.

The rest of the paper is organized as follows.
Section~II provides 
background and problem statements.
Section~III develops the differential privacy
mechanisms for the Laplacian spectrum.
Next, Section~IV compares the scaling of edge and node differential privacy and as a result we shift our attention to edge privacy exclusively. Then, we use the output of the edge mechanism to bound other
graph properties in Section~V. 
Section~VI provides
guidelines and examples and Section~VII concludes. 

\textbf{Notation} 
We use~$\R$ and~$\N$ to denote the real and natural numbers, respectively. 
We use~$|S|$ to denote the cardinality
of a finite set~$S$, and we use~$S_1 \Delta S_2 = (S_1 \backslash S_2) \cup (S_2 \backslash S_1)$
to denote the symmetric difference of two sets. 
For~$n \in \N$, we use~$\mathcal{G}_n$ to denote the set of graphs on~$n$ nodes. 

%% file: background.tex
\section{Preliminaries and Problem Statement}

\subsection{Graph Theory Background}
We consider an undirected, unweighted graph~$G = (V,E)$ defined over a set of nodes $V = \{1,\dots, n \}$
with edge set $E\subset V\times V$. The pair $(i,j)$ belongs to $E$ if nodes $i$ and $j$ share an edge, and $(i,j) \notin E$ otherwise. 
We let $d_i = \vert \{ j\in V \mid (i,j) \in E\}\vert$ denote the degree of node $i\in V$. The degree matrix $D(G)\in \mathbb{R}^{n\times n}$ is the diagonal matrix $D(G)=\operatorname{diag}\big(d_{1}, \ldots, d_{n}\big)$. The adjacency matrix of~$G$ is
\begin{equation}
    (H(G))_{i j}=\begin{cases}
1 & (i, j) \in E \\
0 & \text { otherwise }
\end{cases}.
\end{equation}
We denote the Laplacian of graph $G$ by~${L(G)=D(G)-H(G)}$, which we simply write as $L$ when the associated graph is clear from  context. 

Let the eigenvalues of~$L$ be ordered
according to~${\lambda_1 \leq \lambda_2 \leq \cdots \leq \lambda_n}$. 
The matrix~$L$ is symmetric and positive semidefinite, and thus~$\lambda_i \geq 0$ for all~$i$. 
All graphs~$G$ have~$\lambda_1 = 0$, and a seminal result shows that~$\lambda_2 > 0$
if and only if~$G$ is connected \cite{fiedler1975property}. 
Thus,~$\lambda_2$ is often called the \emph{algebraic connectivity} of a graph. Throughout this paper,
we consider connected graphs with $n\geq3$. 

The value of~$\lambda_2$ specifically encodes a great deal of
information about~$G$: its value is non-decreasing in the number
of edges in~$G$, and algebraic connectivity is closely related to graph diameter and various other algebraic
properties of graphs~\cite{abreu07}. The value of~$\lambda_2$ also characterizes the performance
of consensus algorithms. Specifically, worst-case disagreement in a consensus protocol decays proportionally to~$e^{-\lambda_2 t}$ \cite{Mesbahi2010}. 
Thus, we will privatize the full spectrum of~$L$ and pay special attention to~$\lambda_2$ as we do so.


\subsection{Privacy Background}
We follow the differential privacy definition in~\cite{dwork_algorithmic_2013}. Differential privacy is enforced by a \emph{mechanism}, which is a randomized map. 
Given ``similar" inputs, a differential privacy mechanism produces outputs that are approximately indistinguishable from each other. 
Formally, a mechanism must obfuscate differences between inputs that are 
\emph{adjacent}\footnote{The word ``adjacency'' appears in two forms in this paper: for the adjacency matrix~$H$ above,
and for the adjacency relation used by differential privacy. The adjacency matrix appears only in this section and only
to define the graph Laplacian, and all subsequent uses of ``adjacent'' and ``adjacency'' pertain
to differential privacy (not the adjacency matrix).}. In this work, we analyze two different notions of adjacency for a given graph $G$: an adjacency relation defined with respect to the edges of $G$, $E(G)$, and an adjacency relation defined with respect to the nodes of $G$, $V(G)$. When adjacency is defined with respect to the edge set, we will calibrate our privacy to obfuscate the absence or presence of one or more edges in $G$. When adjacency is defined with respect to the node set, we will obfuscate the absence or presence of a single node. Mathematically, this is done as follows.

\begin{subtheorem}{definition}\label{dfn:both_adj}
\begin{definition}[\emph{Edge Adjacency relation}]
    \label{dfn:edge_adjacency}
Let~$A \in \mathbb{N}$ be given, and
    fix a number of nodes~$n \in \mathbb{N}$. Two graphs~$G, G' \in \mathcal{G}_n$ are adjacent
    if they differ by $A$ edges. We express this mathematically via
    \begin{equation}
    \adj_{e,A}(G, G') = \begin{cases} 1 & |E(G) \Delta E(G')| \leq A \\
                           0 & \textnormal{otherwise}
               \end{cases}.
               \tag*{$\lozenge$}
    \end{equation}
    
\end{definition}

\begin{definition}[\emph{Node Adjacency relation}]
\label{dfn:node_adjacency}
Fix $n\in\mathbb{N}$. Two graphs, $G, G' \in \mathcal{G}_n$ are adjacent 
if they differ by one node with the corresponding edges added or deleted. We express this mathematically via
\begin{equation}
\text{Adj}_n(G,G')=\begin{cases}
1 & |V(G)\Delta V(G')|\leq1\\
0 & \text{otherwise}
\end{cases}.
\tag*{$\lozenge$}
\end{equation}
\end{definition}
\end{subtheorem}

In Definition \ref{dfn:edge_adjacency}, $A$ is the number of edges whose absence or presence must be concealed by privacy, while Definition~\ref{dfn:node_adjacency} specifies that the absence or presence of a single node must be concealed by privacy. In Section~\ref{sec:node_privacy}, we show that a mechanism that obfuscates the
absence or
presence of only a single node is not practical for large networks and engineering examples, 
and therefore we do not consider obfuscating the presence of arbitrary numbers of nodes.

Next, we briefly review differential privacy; see~\cite{dwork_algorithmic_2013} for a complete exposition.
A privacy mechanism $\mathcal{M}$ for a function~$f$ can be obtained by 
first computing the function $f$ on a given input $x$, and then adding noise to~$f(x)$. 
The distribution of noise depends on the sensitivity of the function $f$ to changes in its input, described below.
It is the role of a mechanism to approximate functions
of sensitive data with private responses, and we next state this formally. The guarantees of privacy are defined with respect to the adjacency relation. 
Since we consider two notions of adjacency, we define 
two types of privacy: 
(i) edge differential privacy using the standard definition of differential privacy equipped with the edge adjacency relation, $\adj_{e,A},$ appearing in Defintion~\ref{dfn:edge_adjacency}, 
and (ii) node differential privacy using the standard definition of differential privacy equipped with $\adj_n$ in Definition~\ref{dfn:node_adjacency}.

\begin{subtheorem}{definition}\label{dfn:both_dp}
\begin{definition}[\emph{Edge differential privacy; \cite{dwork_algorithmic_2013}}]\label{dfn:dp_edge}
    Let $\epsilon > 0$, $\delta \in [0,1)$ be given, use $\adj_{e,A}$ from Definition~\ref{dfn:edge_adjacency}, and fix a probability space $(\Omega, \mathcal{F}, \mathbb{P})$. 
    Then a mechanism $\mathcal {M}: \Omega \times \mathcal{G}_n \rightarrow \mathbb{R}$ is $(\epsilon , \delta)$-differentially private if,
    for all adjacent graphs $G,G' \in \mathcal{G}_n$, 
    \begin{equation}
        \mathbb{P}\big[\mathcal{M}(G) \in S\big] \leq \exp (\epsilon) \cdot \mathbb{P}\big[\mathcal{M}\left(G'\right) \in S\big] + \delta
    \end{equation}
    for all sets $S$ in the Borel $\sigma$-algebra over $\mathbb{R}$. \hfill $\lozenge$
\end{definition}

\begin{definition}[\emph{Node differential privacy; \cite{dwork_algorithmic_2013}}]\label{dfn:dp_node}
    Let $\epsilon > 0$, $\delta \in [0,1)$ be given, use $\adj_n$ from Definition~\ref{dfn:node_adjacency}, and fix a probability space $(\Omega, \mathcal{F}, \mathbb{P})$. 
    Then a mechanism $\mathcal {M}: \Omega \times \mathcal{G}_n \rightarrow \mathbb{R}$ is $(\epsilon , \delta)$-differentially private if,
    for all adjacent graphs $G,G' \in \mathcal{G}_n$, 
    \begin{equation}
        \mathbb{P}\big[\mathcal{M}(G) \in S\big] \leq \exp (\epsilon) \cdot \mathbb{P}\big[\mathcal{M}\left(G'\right) \in S\big] + \delta
    \end{equation}
    for all sets $S$ in the Borel $\sigma$-algebra over $\mathbb{R}$. \hfill $\lozenge$
\end{definition}
\end{subtheorem}

The value of~$\epsilon$ controls the amount of information shared, and typical values range from~$0.1$ to~$\log 3$~\cite{dwork_algorithmic_2013}.
The value of~$\delta$ can be regarded as the probability that more information is shared than~$\epsilon$ should allow, and typical
values range from~$0$ to~$0.05$. Smaller values of both imply stronger privacy. 
Given~$\epsilon$ and~$\delta$, a privacy mechanism must enforce Definition~\ref{dfn:dp_edge} or~\ref{dfn:dp_node} for all
graphs adjacent in the sense of Definition~\ref{dfn:edge_adjacency} or~\ref{dfn:node_adjacency}, respectively. 

We next define the sensitivity of~$\lambda_i$, which will be used later to calibrate the variance of privacy noise.  
With a slight abuse of notation, we treat~$\lambda_i$ as a function~$\lambda_i: \mathcal{G}_n \to \mathbb{R},$ and we will develop differential privacy mechanisms to approximate each $\lambda_i$. 
The sensitivity will depend on which adjacency relation is used, and this is made explicit in the following definitions.

\begin{subtheorem}{definition}\label{dfn:both_sens}
\begin{definition}[\emph{Edge Sensitivity}]\label{dfn:edge_sens}
    The edge sensitivity of~$\lambda_i$ is the greatest difference
    between its values on Laplacians of graphs that are adjacent with respect to $\adj_{e,A}$ in Defintion~\ref{dfn:edge_adjacency}. Formally, for a fixed $A,$ the edge sensitivity of $\lambda_i$ is given as
    \begin{equation}
    \Delta \lambda_{i,e} = \max_{\substack{G, G' \in \mathcal{G}_n \\ \adj_{e,A}(G, G') = 1}} \big|\lambda_i(L) - \lambda_i(L')\big|,
    \end{equation}
    where $L$ and $L'$ are the Laplacians of $G$ and $G'$. \hfill $\lozenge$
\end{definition}

\begin{definition}[\emph{Node Sensitivity}]\label{dfn:node_sens}
    The node sensitivity of~$\lambda_i$ is the greatest difference
    between its values on Laplacians of graphs that are adjacent with respect to $\adj_{n}$ in Defintion~\ref{dfn:node_adjacency}. Formally, the node sensitivity of $\lambda_i$ is given as
    \begin{equation}
    \Delta \lambda_{i,n} = \max_{\substack{G, G' \in \mathcal{G}_n \\ \adj_n(G, G') = 1}} \big|\lambda_i(L) - \lambda_i(L')\big|,
    \end{equation}
    where $L$ and $L'$ are the Laplacians of $G$ and $G'$. \hfill $\lozenge$
\end{definition}
\end{subtheorem}
Noise is added by a mechanism, which is a randomized map used to implement differential privacy. 
The Laplace mechanism is widely used, and it adds noise from a Laplace distribution to sensitive data (or functions thereof). 
The standard Laplace mechanism has support on all of~$\mathbb{R}.$ For graphs on~$n$ nodes, $\lambda_i\in[0, n]$~for all $i.$ 
To generate a private output, one can add Laplace noise and then project the result onto~$[0, n]$ (which is differentially private because the projection is post-processing), though
similar approaches have been shown to produce highly inaccurate private data~\cite{gohari21}. 
Instead, we use the bounded Laplace mechanism in~\cite{holohan2018bounded}. 
We state it in a form amenable to use with~$\lambda_i$. 

\begin{definition}    \label{dfn:bounded_laplace}
 Let $b >0$ and let $D = [0, n]$. 
	Then the bounded Laplace mechanism $W_{\lambda_i}:\Omega \rightarrow D$, for each
	$\lambda_i\in D$, is given by its probability density function $f_{W_{\lambda_i}}$ as
	\begin{equation}
		f_{W_{\lambda_i}}(x)
		=\begin{cases}
			0 & \text{if }x\notin D\\
			\frac{1}{C(\lambda_i,b)}\frac{1}{2b}e^{-\frac{|x-\lambda_i|}{b}} & \text{if }x\in D
		\end{cases},
	\end{equation}
	\noindent where $C(\lambda_i,b)=\int_{D}\frac{1}{2b}e^{-\frac{|x-\lambda_i|}{b}}dx.$ \hfill$\lozenge$
\end{definition}
\subsection{Problem Statements}
We now give formal problem statements. The first two pertain to the development of privacy mechanisms. 
\begin{problem}
\label{prb:prob1}
Develop a mechanism to provide $(\epsilon, \delta)$-edge differential privacy 
in the sense of Definition~\ref{dfn:dp_edge} 
for the spectrum of the graph Laplacian $L(G)$ of a graph $G.$
\end{problem}

\begin{problem}
\label{prb:prob2}
Develop a mechanism to provide $(\epsilon, \delta)$-node differential privacy 
in the sense of Definition~\ref{dfn:dp_node} 
for the algebraic connectivity of a graph $G$.
\end{problem}

We note that Problem~\ref{prb:prob2} considers the algebraic connectivity specifically because
that will be used to show the poor scaling of node privacy for the full Laplacian spectrum. 
Comparisons of the two mechanisms are the subject of the next problem. 

\begin{problem}
\label{prb:prob3}
Given a graph $G$ on $n$ nodes and two privacy mechanisms, $\mathcal{M}_n$ and $\mathcal{M}_e$, that provide $(\epsilon,\delta)$-node privacy and $(\epsilon, \delta)$-edge privacy for the spectrum of $L(G),$ respectively, 
analyze how the variances of the two mechanisms scale with respect to the size of the network $n.$
\end{problem}

The final two problem statements pertain to the accuracy of graph properties when bounded
using private spectra. 

\begin{problem}
\label{prb:prob5}
Given a private algebraic connectivity, develop
bounds on the expectation of the graph diameter and mean distance between nodes in the graph.
\end{problem}

\begin{problem}
\label{prb:prob4}
Given private values of the Laplacian spectrum, provide examples to numerically quantify the accuracy of using these private values to estimate the trace of the Laplacian, Kemeny's constant, and Cheeger's inequality.
\end{problem}

%% file: privacy_mechanisms.tex
\section{Privacy Mechanisms}
\label{sec:privacy_mechanisms}
In this section, we solve Problems~\ref{prb:prob1} and~\ref{prb:prob2}. Specifically, we develop two mechanisms to provide $(\epsilon,\delta)-$differential privacy to eigenvalues of a graph Laplacian $L.$  In Section~\ref{sec:edge_privacy}, we use edge differential privacy to privatize each of the Laplacian eigenvalues, $\lambda_i$ for $i\in[n].$ 
Then in Section~\ref{sec:node_privacy}, we use node differential privacy to privatize $\lambda_2.$ In both subsections we first bound the sensitivity appearing in Definition~\ref{dfn:both_sens} and then use these sensitivity bounds to develop the privacy mechanisms.

\subsection{Edge Privacy}
\label{sec:edge_privacy}
We now design a mechanism to implement $(\epsilon,\delta)-$edge differential privacy. We first bound the sensitivity $\Delta \lambda_{i,e}$ appearing in Definition~\ref{dfn:edge_sens}.  

\begin{lemma}[Edge sensitivity bound]
\label{lem:edge_sensitivity}
Fix an adjacency parameter~$A \in \N$. Then for the edge sensitivity~$\Delta \lambda_{i,e}$ in Definition~\ref{dfn:edge_sens}, we have 
\begin{equation}
\Delta\lambda_{i,e} \leq 2A
\end{equation}
for $i \in \{1, \dots, n\}.$
\end{lemma}
\emph{Proof:} See Appendix~A.\hfill$\blacksquare$

Next, we establish an algebraic relation for $b_e$, which lets the bounded Laplace mechanism satisfy the theoretical guarantees  of $(\epsilon, \delta)$-edge differential privacy in Definition~\ref{dfn:dp_edge}.
\begin{theorem}
\label{thm:mechanim_variance_ub_edge}
Let~$\epsilon > 0$ and~$\delta \in (0, 1)$ be given.
Fix~$n \in \N$ and consider graphs in~$\mathcal{G}_n$. 
Then for the bounded Laplace mechanism~$W_{\lambda_i}$ in 
Definition~\ref{dfn:bounded_laplace}, choosing $b_e$ according to 
\begin{equation} \label{eq:bdef}
    	b_e \geq \frac{2 A}{\epsilon-\log \left(\frac{2-e^{-\frac{2 A}{b_e}}-e^{-\frac{n-2 A}{b_e}}}{1-e^{-\frac{n}{b_e}}}\right)-\log (1-\delta)}
\end{equation}
satisfies $(\epsilon , \delta)$-edge differentially privacy with respect to $\adj_{e,A}$ as defined in Definition~\ref{dfn:dp_edge}.
\end{theorem}
\emph{Proof:} By~\cite[Theorem 3.5]{holohan2018bounded}, the bounded Laplace
mechanism provides $(\epsilon,\delta)-$differential privacy if 
\begin{equation}
	b_e\ge\frac{\Delta \lambda_{i,e}}{\epsilon-\log\Delta C(b_e)-\log(1-\delta)},
\end{equation}
where, given that $\lambda_i\in [0,n]$, $ \Delta C(b_e)$ is defined as
\begin{equation}\label{eq:delta_C_b_edge}
\Delta C(b_e) :=\frac{C(\Delta \lambda_{i,e},b_e)}{C(0,b_e)},
\end{equation}
where~$C$ is from Definition~\ref{dfn:bounded_laplace}. 
Next, we find

\begin{align}
    C(\lambda_{i},b_e) &= \int_0^n \frac{1}{2b_e}e^{-\frac{|x-\lambda_i|}{b_e}}dx \\
    & =\frac{1}{2b_e}\int_0^{\lambda_i}e^{\frac{x-\lambda_i}{b_e}}dx + \frac{1}{2b_e}\int_{\lambda_i}^n e^{-\frac{x-\lambda_i}{b_e}}dx\\
    &=1-\frac{1}{2}\left(e^{-\frac{\lambda_{i}}{b_e}}+e^{-\frac{n-\lambda_{i}}{b_e}}\right) \label{eq:cl2_edge}. 
\end{align}
Using~\eqref{eq:cl2_edge} to compute~$C(\Delta\lambda_{i,e}, b_e)$ and~$C(0, b_e)$ in~\eqref{eq:delta_C_b_edge} gives 
\begin{equation}
	\Delta C(b_e) =\frac{1-\frac{1}{2}\left(e^{-\frac{\Delta \lambda_{i,e}}{b_e}}+e^{-\frac{n-\Delta \lambda_{i,e}}{b_e}}\right)}{1-\frac{1}{2}\left(1+e^{-\frac{n}{b_e}}\right)}. 
\end{equation}
Using the sensitivity bound in Lemma~\ref{lem:edge_sensitivity}, we put $\Delta \lambda_{i,e}=2A$, which completes the proof. \hfill $\blacksquare$

Theorem~\ref{thm:mechanim_variance_ub_edge} solves Problem~\ref{prb:prob1}, 
and we now have an $(\epsilon,\delta)-$edge differential privacy mechanism for the spectrum of a graph Laplacian $L$. We now shift our attention to node privacy and the algebraic connectivity, $\lambda_2.$

\subsection{Node Privacy}
\label{sec:node_privacy}
Here we develop an $(\epsilon,\delta)-$node differential privacy mechanism for $\lambda_2.$ We will use the same process as the last subsection: we first bound $\Delta \lambda_{2,n}$ from Definition~\ref{dfn:node_sens} for a graph $G\in\mathcal{G}_n$, then use this sensitivity to find an algebraic relation for the bounded Laplace mechanism to satisfy Definition~\ref{dfn:dp_node}. In Lemma~\ref{lem:edge_sensitivity}, we were able to derive a common bound on the sensitivity of each eigenvalue of $L$ when edge sensitivity is used. There is no common bound when node sensitivity is used. 
In Section~\ref{sec:scaling}, we show that the node privacy scales poorly with the size of the network and will not be usable in most engineering problems. 
Thus, in this section we focus on $\lambda_2$ rather than the entire spectrum, as this is sufficient to illustrate the poor scaling of node
privacy in this context. 

Definition~\ref{dfn:node_adjacency} considers adjacent graphs as graphs that have an additional or absent node from $G.$ Thus, for a $G'$ satisfying $\adj_n(G,G')=1,$ it is possible that  ${G'\in\mathcal{G}^{n-1}}$ or ${G'\in\mathcal{G}^{n+1}}.$ Because of this, we require $n\geq 3$ and the two cases will be handled separately in our analysis. 
We have the following result.

\begin{lemma}
\label{lem:node_sensitivity}
Fix~$n \in \N$ and consider graphs in~$\mathcal{G}_n$. Then the node sensitivity of $\lambda_2$ in Definition~\ref{dfn:node_sens} is bounded as
\[
\Delta\lambda_{2,n}\leq n-1.
\]
\end{lemma}
\emph{Proof:} See Appendix~B.\hfill$\blacksquare$

With this sensitivity bound, we now establish an algebraic relation for $b_n$, which lets the bounded Laplace mechanism satisfy the theoretical guarantees  of $(\epsilon, \delta)$-node 
differential privacy in Definition~\ref{dfn:dp_node}.

\begin{theorem}
\label{thm:mechanim_variance_ub_node}
Let~$\epsilon > 0$ and~$\delta \in (0, 1)$ be given. 
Fix~$n \in \N$ and consider graphs in~$\mathcal{G}_n$. 
Then for the bounded Laplace mechanism~$W_{\lambda_2}$ in 
Definition~\ref{dfn:bounded_laplace}, choosing $b_n$ according to 
\begin{equation} \label{eq:bdef}
    	b_n\geq\frac{n-1}{\epsilon-\log\left(\frac{2-e^{-\frac{n-1}{b_n}}-e^{-\frac{1}{b_n}}}{1-e^{-\frac{n}{b_n}}}\right)-\log(1-\delta)}
\end{equation}
satisfies $(\epsilon , \delta)$-node differential privacy with respect to $\text{Adj}_n$ from Definition~\ref{dfn:node_adjacency}.
\end{theorem}
\emph{Proof:} By~\cite[Theorem 3.5]{holohan2018bounded}, the bounded Laplace
mechanism satisfies differential privacy if 
\begin{equation}
	b_n\geq\frac{\Delta\lambda_{2,n}}{\epsilon-\log\Delta C(b_n)-\log(1-\delta)},
\end{equation}
where, given that $\lambda_2\in [0,n]$, $ \Delta C(b_n)$ is defined as
\begin{equation}\label{eq:delta_C_b}
\Delta C(b_n) :=\frac{C(\Delta \lambda_{2,n},b_n)}{C(0,b_n)},
\end{equation}
where~$C$ is from Definition~\ref{dfn:bounded_laplace}. 
Next, we find
\begin{equation} \label{eq:cl2}
	C(\lambda_{2},b_n)=1-\frac{1}{2}\left(e^{-\frac{\lambda_{2}}{b_n}}+e^{-\frac{n-\lambda_{2}}{b_n}}\right).
\end{equation}
Using~\eqref{eq:cl2} to compute~$C(\Delta\lambda_{2, n}, b_n)$ and~$C(0, b_n)$ in~\eqref{eq:delta_C_b} gives 
\begin{equation}
	\Delta C(b_n) =\frac{1-\frac{1}{2}\left(e^{-\frac{\Delta \lambda_{2,n}}{b_n}}+e^{-\frac{n-\Delta \lambda_{2,n}}{b_n}}\right)}{1-\frac{1}{2}\left(1+e^{-\frac{n}{b_n}}\right)}. 
\end{equation}
Using the sensitivity bound in Lemma~\ref{lem:node_sensitivity}, we set $\Delta \lambda_{2,n}=n-1$, which completes the proof. \hfill $\blacksquare$

Theorem~\ref{thm:mechanim_variance_ub_node} solves Problem~\ref{prb:prob2} and gives an $(\epsilon,\delta)-$node differential privacy mechanism for the algebraic connectivity, $\lambda_2$, of the graph Laplacian $L$. The algebraic relationships given in Theorems~\ref{thm:mechanim_variance_ub_edge} and~\ref{thm:mechanim_variance_ub_node} are defined implicitly since $b$ appears on both sides of the expression, and they do not yield an analytical expression for the minimum $b$ required for privacy. In \cite{holohan2018bounded}, the authors provide an algorithm to solve for~$b$ using the bisection method, and we use this in the remainder of this paper. However, the lack of an analytical expression for the required $b$ prevents us from immediately comparing the 
amounts of noise required by the 
two notions of privacy. The next section derives necessary conditions 
for the variances of noise required
for edge and node privacy, which will allow us to compare how the two notions of privacy scale with the size of the network $n$. 


%% file: scaling.tex
\section{Scaling Laws}
\label{sec:scaling}
In this section we will compare the notions of edge and node differential privacy to solve Problem~\ref{prb:prob3}. 
More specifically, we will analyze how the required variance of each privacy notion scales with the size of the network $n.$ Here, we focus on the algebraic connectivity $\lambda_2$ to draw accurate comparisons between edge and node privacy. However, the edge privacy results can immediately 
be applied to the rest of the Laplacian's spectrum and the scaling trends found here persist for each value of the Laplacian spectrum.
\begin{figure}
    \centering
    \includegraphics[width=.5\textwidth,draft=false]{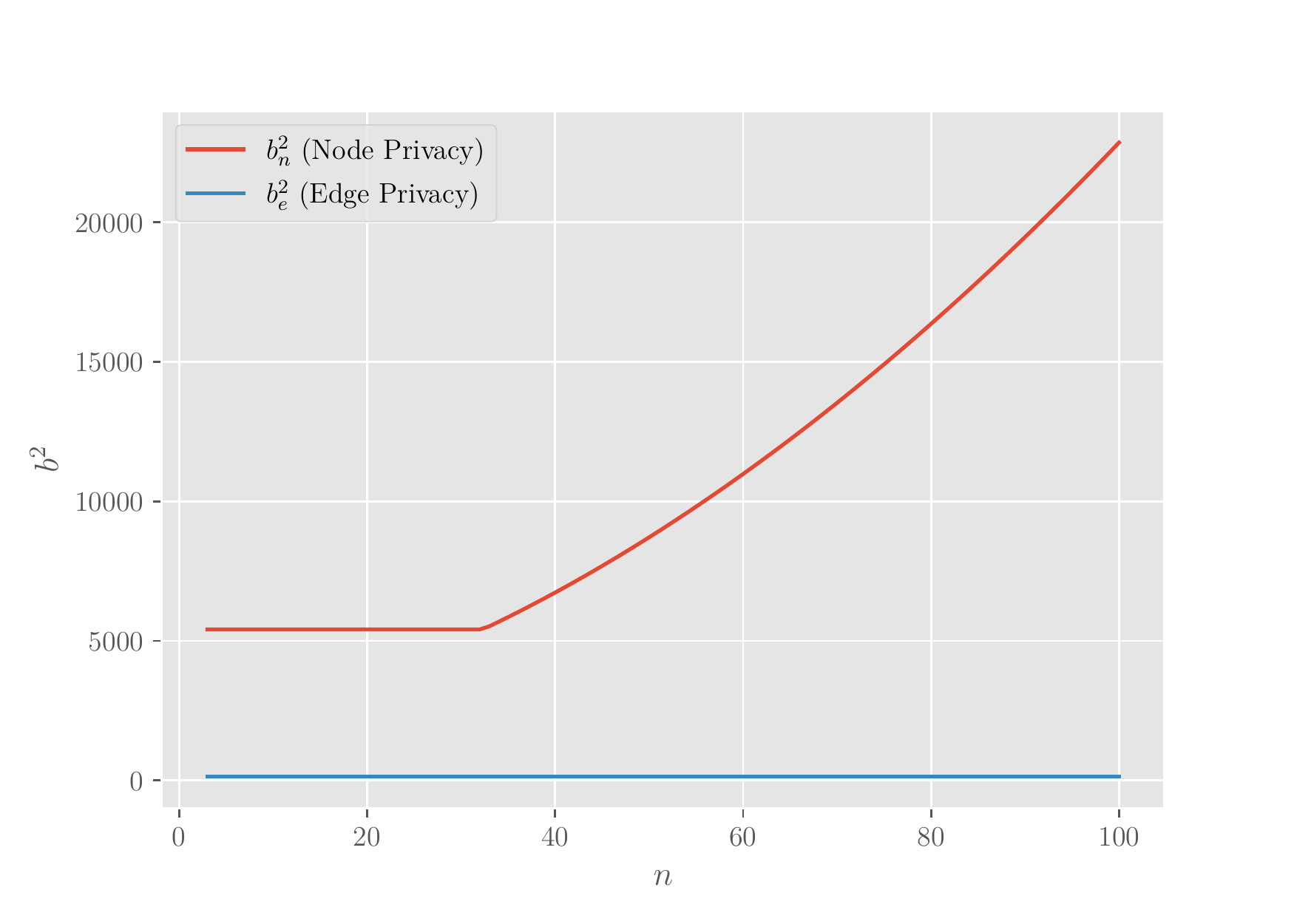}
    \caption{Fix $\epsilon=0.4,\delta=0.05,A=2,$ and $\lambda_2=2.5$. 
    We set~$b_n$ equal to its lower bound in Corollary~\ref{cor:weaker_node_privacy} and~$b_e$ equal to its lower bound in Corollary~\ref{cor:weaker_edge_privacy}.
    The variances of the Laplace mechanisms are proportional to~$b_e^2$ and~$b_n^2$, and we plot~$b_e^2$ and~$b_n^2$ here for~$n=3$ to~$n=100$ nodes. 
    This figure shows that the variance of noise required for edge privacy has no dependence on $n$, while it is necessary
    for the variance of noise for node privacy to grow quadratically in~$n$. 
    }
    \label{fig:scaling_plot}
\end{figure}
To compare the two mechanisms, we fix a graph $G\in\mathcal{G}_n$ and privacy parameters $\epsilon$ and $\delta.$ Then we define an edge and node privacy mechanism to provide $(\epsilon,\delta)-$differential 
privacy with parameters $b_e$ and $b_n,$ respectively. Then we will analyze and compare the required values of $b_e$ and $b_n$ given this $\epsilon$ and $\delta.$

\subsection{Comparison of Mechanisms}

Recall that the requirements for the bounded Laplace mechanism to achieve $(\epsilon,\delta)-$differential privacy appearing in Theorems~\ref{thm:mechanim_variance_ub_edge} and~\ref{thm:mechanim_variance_ub_node} are  defined implicitly in $b_e$ and $b_n$ and the minimal values must be found numerically. To compare the two notions of privacy we find weaker, necessary conditions for $(\epsilon,\delta)-$edge and node differential privacy, which give an analytical expression for the growth of $b$. The following results will show that the required parameter for the bounded Laplace mechanism to achieve $(\epsilon,\delta)-$differential privacy is strictly larger than the parameter required for the standard, unbounded Laplace mechanism from~\cite{dwork_algorithmic_2013} to achieve the same level of privacy. This recovers a general-purpose result of the same kind presented in~\cite[Theorem 3.5]{holohan2018bounded}. The next two corollaries give these necessary conditions for edge and node privacy, respectively.

\begin{corollary}
\label{cor:weaker_node_privacy}
Fix a graph $G\in\mathcal{G}^n,$ $\epsilon>0,$ and $\delta\in(0,1).$ Let $W^n_{\lambda_2}$ be a bounded Laplace mechanism with parameter $b_n.$ Then 
\begin{equation} 
    b_n>\frac{n-1}{\epsilon-\log(1-\delta)}
\end{equation}
is a necessary condition for $W^n_{\lambda_2}$ to provide $(\epsilon,\delta)-$node differential privacy.
\end{corollary}
\emph{Proof:}
See Appendix~C.\hfill$\blacksquare$

\begin{corollary}
\label{cor:weaker_edge_privacy}
Fix a graph $G\in\mathcal{G}^n,$ $\epsilon>0,$ and $\delta\in(0,1).$ Let $W^e_{\lambda_2}$ be a bounded Laplace mechanism with parameter $b_e.$ Then 
\begin{equation} 
    b_e>\frac{2A}{\epsilon-\log(1-\delta)}
\end{equation}
is a necessary condition for $W^e_{\lambda_2}$ to provide $(\epsilon,\delta)-$edge differential privacy.
\end{corollary}
\emph{Proof:}
See Appendix~D.\hfill$\blacksquare$

\begin{remark}
In Corollary~\ref{cor:weaker_node_privacy}, the necessary condition on~$b_n$ for $(\epsilon,\delta)$-node differential privacy scales linearly with $n$. 
A standard Laplace distribution with parameter $b$ has variance $2{b}^2.$ This means that as the size of the network $n$ grows, the variance required 
for $(\epsilon,\delta)-$node differential privacy grows quadratically in $n.$ Simultaneously, in Corollary~\ref{cor:weaker_edge_privacy}, $b_e$ has no dependence 
on the size of the network. Thus,~$b_e$ and the variance of privacy noise needed for edge differential privacy remain constant as a network grows. 

\end{remark}

Corollaries~\ref{cor:weaker_node_privacy} and~\ref{cor:weaker_edge_privacy} allow us to make further comparisons between edge and node privacy. For example, fix an $\epsilon,\delta,$ and $n.$ For the edge privacy mechanism $W_{\lambda_2}^e$ to have larger variance than the node privacy mechanism $W_{\lambda_2}^n$ according to the necessary conditions, we would need $A\geq\frac{n-1}{2}.$ This implies that the noise required to obfuscate the presence of a single node is proportional to the noise required to obfuscate the presence of $\frac{n-1}{2}$ edges. In applications where many connections need to obfuscated, specifically on the order of half the connections in the network, node privacy can provide similar accuracy to edge privacy. However, in applications where less than half the connections in the network need to be obfuscated, edge privacy can provide greater accuracy than node privacy.

The aforementioned scaling laws are further illustrated by numerical results shown in Figure~\ref{fig:scaling_plot}. These results show that the minimal $b_n$
required for node privacy grows quickly, while~$b_e$ remains constant. 

One of the appealing features of differential privacy is that it provides a means to share private information that can still can still be useful.  However, in most applications, variance on the order of $n^2$ will render the private information useless. Thus, we focus on edge privacy for the rest of the paper.

\subsection{Accuracy of Edge Privacy}


The edge privacy mechanism has the following accuracy.

\begin{theorem}\label{thm:edge_privacy_acc}
For a fixed $G\in\mathcal{G}_n$, $\epsilon > 0$, $\delta \in (0, 1)$, and $A \in \mathbb{N}$, the accuracy of a private eigenvalue $\tilde\lambda_i$ generated using the bounded Laplace mechanism with parameter $b_e$ is given by
\begin{equation}
  E\left[\tilde{\lambda}_{i}-\lambda_{i}\right] = 
	\frac{1}{2C(\lambda_{i},b_e)}\left(2\lambda_{i}+b_e e^{-\frac{\lambda_{i}}{b_e}}-(n+b_e)e^{-\frac{n-\lambda_{i}}{b_e}}\right)-\lambda_{i}.
\end{equation}
\end{theorem}
\emph{Proof:} See Appendix~E\hfill$\blacksquare$

Theorem~\ref{thm:edge_privacy_acc} provides an analytical expression for the accuracy of the edge privacy mechanism. Since differential privacy is immune to post-processing, we can use private spectra to estimate other graph properties without harming privacy guarantees. The rest of the paper focuses on estimating othering graph properties using the edge privacy mechanism. Specifically, in Section~\ref{sec:bo_section} we develop statistical bounds on other graph properties given a private $\lambda_2,$ and in Section~\ref{sec:examples} we provide a series of examples that demonstrate the accuracy of the edge privacy mechanism and illustrate how these private values of the Laplacian spectrum can be used to estimate other graph properties.
    
    

%% file: other_properties.tex
\section{Bounding Other Graph Properties}
\label{sec:bo_section}
In this section we solve Problem~\ref{prb:prob5}. There exist numerous inequalities relating~$\lambda_2$ to
other quantitative graph properties~\cite{abreu07,Mesbahi2010}, 
and one can therefore expect that the private~$\lambda_2$ will be used to estimate
other quantitative characteristics of graphs. 
To illustrate the utility of doing so, 
in this section we bound the
graph diameter~$d$ and mean distance~$\rho$
in terms of the private value~$\tilde{\lambda}_2$. 

\subsection{Analytical bounds}
Both~$d$ and~$\rho$ measure graph size and provide insight into 
how easily information can be transferred across a network~\cite{PALDINO2017201}. 
We estimate each one in terms of the private~$\lambda_2$ and bound the error
induced in these estimates by privacy. 
These bounds represent the types of calculations one can do with~$\tilde{\lambda}_2$, and
similar bounds can be easily derived, e.g., on minimal/maximal degree, 
edge connectivity, etc., because their bounds are proportional to~$\lambda_2$~\cite{fiedler_algebraic_1973}.

We first recall bounds from the literature. 
\begin{lemma}[Diameter and Mean Distance Bounds\cite{Mohar1991eigenvalues}]\label{thm:d_rho_bounds}
    For an undirected, unweighted graph $G$ on $n$ nodes, define 
    \begin{align*}
        &\overline d(\lambda_2,\alpha) =  \left(2\sqrt{\frac{\lambda_n}{\lambda_2}}\sqrt{\frac{\alpha^2-1}{4\alpha}}+2\right)\left(\log_\alpha\frac{n}{2}\right) \\
        &\overline\rho(\lambda_2,\alpha) =  \left(\!\sqrt{\frac{\lambda_n}{\lambda_2}}\sqrt{\frac{\alpha^2-1}{4\alpha}}\!+\!1\!\right)\left(\!\frac{n}{n-1}\!\right)\left(\frac{1}{2} \!+\! \log_\alpha\frac{n}{2}\right).
    \end{align*}
    Then for any fixed~$\lambda_2>0$ and any~$\alpha>1$,
    the diameter~$d$ and mean distance~$\rho$ of the graph~$G$ are bounded via
    \begin{align}
        &\underline{d}(\lambda_2)=\frac{4}{n\lambda_2}\leq d \leq  \overline d(\lambda_2,\alpha) \\
        &\underline{\rho}(\lambda_2)=\frac{2}{(n-1)\lambda_2} + \frac{n-2}{2(n-1)} \leq \rho \leq \overline\rho(\lambda_2,\alpha). 
    \end{align}
    The least upper bounds can be derived by finding values of $\alpha_d$ and $\alpha_\rho$ which minimize $\overline d(\lambda_2,\alpha)$ and $\overline\rho(\lambda_2,\alpha)$, respectively.\hfill $\blacksquare$
\end{lemma} 

A list of $\alpha_d$ and $\alpha_\rho$ values can be found in Table 1 in \cite{Mohar1991eigenvalues}. To quantify the impacts of 
using the private~$\lambda_2$ in these bounds, we next bound the expectations of the private forms of~$d$ and~$\rho$. 
These bounds use the upper incomplete gamma function~$\Gamma(\cdot,\cdot)$ and the imaginary error function~$\textrm{erfi}(\cdot)$,
defined as
\begin{equation}
\Gamma(s,x) = \int_x^\infty t^{s-1}e^{-t}dt \,\,\,\textnormal{ and }\,\,\, \textrm{erfi}(x) = \frac{2}{\sqrt{\pi}}\int_0^x e^{t^2}dt.
\end{equation}

Using the private~$\lambda_2$, expectation bounds are as follows. 
\begin{theorem}[\emph{Expectation bounds for~$d$ and~$\rho$; Solution to Problem~4}]\label{thm:expectation_bounds_d_rho}
    For any $\lambda_2> 0$, denote its private value by~$\tilde{\lambda}_2$. 
    Let~$\tilde{d}$ and~$\tilde{\rho}$ denote the estimates of the diameter and mean distance, respectively, when computed with~$\tilde{\lambda}_2$. 
    Then the expectations~$E[\tilde{d}]$ and~$E[\tilde{\rho}]$, obey
    \begin{align}
        &\frac{4}{nE[\tilde{\lambda}_2]}\leq E[\tilde{d}] \leq E[\overline{d}(\tilde{\lambda}_2,\alpha_d)] \qquad \textnormal{ and }\\
        &\frac{2}{(n-1)E[\tilde{\lambda}_2]} + \frac{n-2}{2(n-1)} \leq E[\tilde{\rho}]\leq E[\overline{\rho}(\tilde{\lambda}_2,\alpha_\rho)],
    \end{align}
    where
    \begin{align*}
        &E[\overline{d}(\tilde{\lambda}_2,\alpha_d)] \!=\!\left[2\sqrt{\frac{\lambda_n(\alpha_d^2-1)}{4\alpha_d}}E\!\left[\sqrt{\frac{1}{\tilde{\lambda}_2}}\right] \!+\! 2\right]\!\!\left[\log_{\alpha_d}\frac{n}{2}\right] \\
        &\begin{multlined}E[\overline{\rho}(\tilde{\lambda}_2,\alpha_\rho)]=\left[\sqrt{\frac{\lambda_n(\alpha_\rho^2-1)}{4\alpha_\rho}}E\left[\frac{1}{\sqrt{\tilde{\lambda}_2}}\right]+1\right]\\
        \cdot\left[\frac{n}{n-1}\right]\cdot\left[\frac{1}{2}+\log_{\alpha_\rho}\frac{n}{2}\right]\end{multlined}.
    \end{align*}
    We can compute the expectation terms with~$\tilde{\lambda}_2$ via
    \begin{align*}
        &\begin{multlined}
        E\left[\frac{1}{\sqrt{\tilde{\lambda}_2}}\right] = \frac{1}{C(\lambda_2,b_e)}\frac{1}{2b_e}\left(\sqrt{\pi}\sqrt{b_e}e^{-\frac{\lambda_2}{b_e}}\left(\textrm{erfi}\left(\sqrt{\frac{\lambda_2}{b_e}}\right)\right)\right.\\
        \left.+\sqrt{b_e}e^{\frac{\lambda_2}{b_e}}\left(\Gamma\left(\frac{1}{2},\frac{n}{b_e}\right)-\Gamma\left(\frac{1}{2},\frac{\lambda_2}{b_e}\right)\right)\right)
        \end{multlined} \\
        &\begin{multlined}
        E[\tilde{\lambda}_2] = \frac{1}{2C(\lambda_2,b_e)}\left(2\lambda_2+b_e e^{-\frac{\lambda_2}{b_e}}-b_e e^{-\frac{n-\lambda_2}{b_e}}-ne^{-\frac{n-\lambda_2}{b_e}}\right),
        \end{multlined}
    \end{align*}
    where~$C$ is from Definition~\ref{dfn:bounded_laplace}. 
\end{theorem}

\emph{Proof:} 
See Appendix~F.
\hfill $\blacksquare$

\begin{remark}\label{rmk:epsilon_bounds_relation}
    A larger~$\epsilon$ gives weaker 
    privacy, and it results in a smaller value of $b_e$ and a distribution
    of privacy noise that is more tightly concentrated about its mean. 
    Thus, a larger~$\epsilon$ implies that the expected value
    $E[\tilde{\lambda}_2]$ is closer to the exact,
    non-private $\lambda_2$, which leads to smaller disagreements in the bounds on 
    the true and expected values of $d$ and $\rho$.     
\end{remark}

\subsection{Simulation results}

\begin{figure}
    \centering
    \includegraphics[width=.45\textwidth,draft=false]{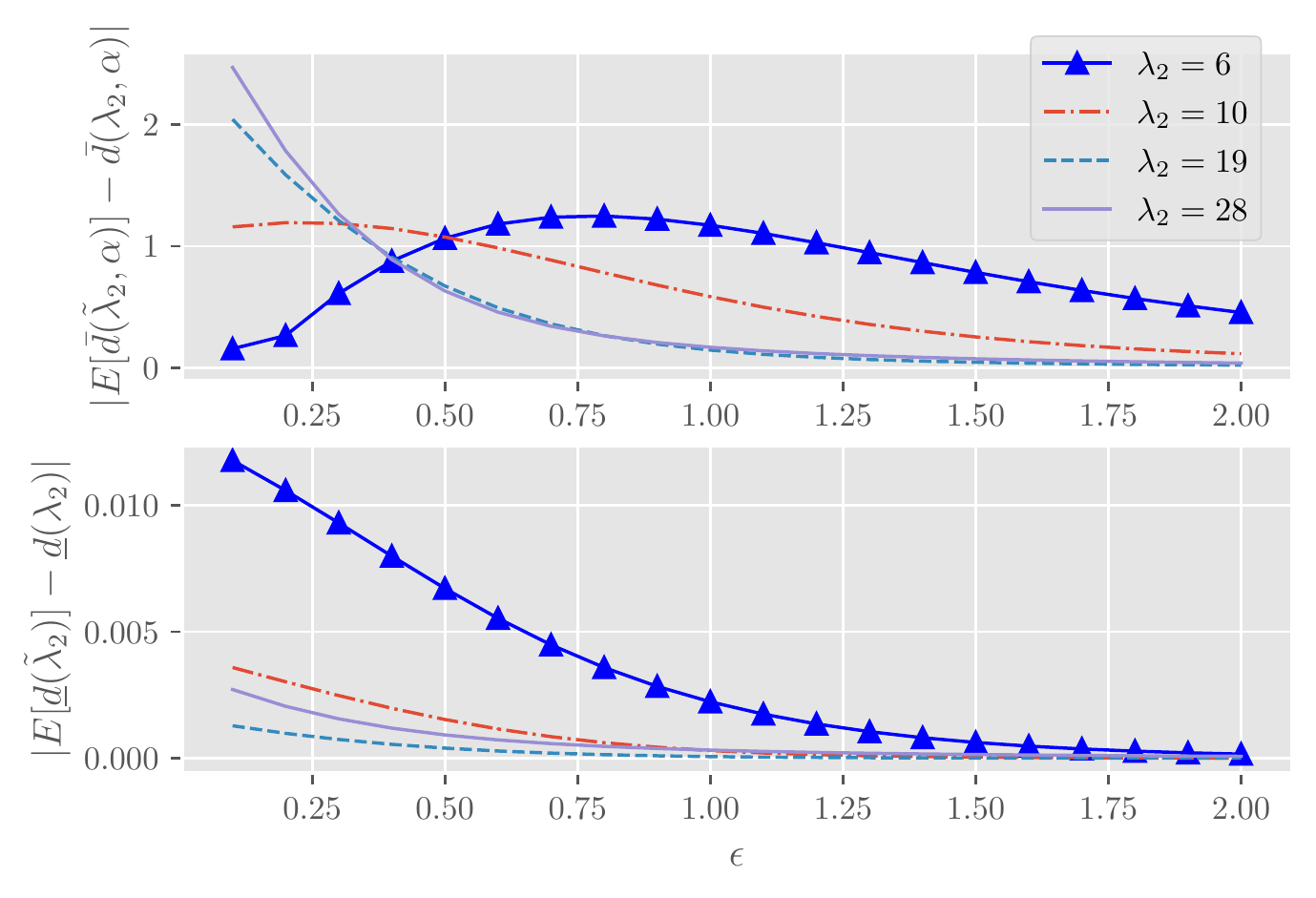}
    \vspace{-5mm}
    \caption{The top plot shows the distance between the exact and expected upper bounds for~$d$. The bottom plot shows the distance between the corresponding lower bounds.}
    \label{fig:bounds_for_d}
\end{figure}

\begin{figure}
    \centering
    \includegraphics[width=.45\textwidth,draft=false]{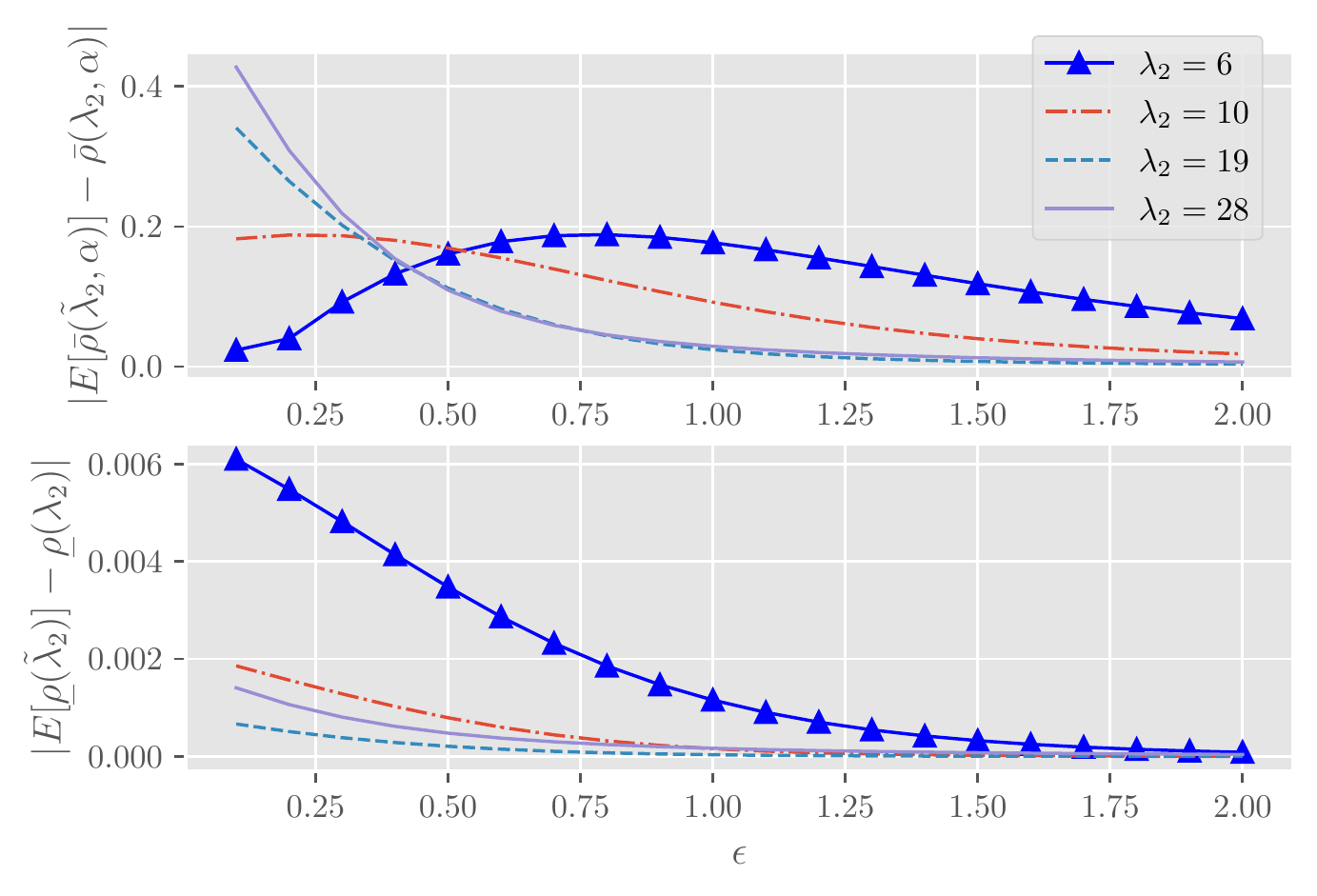}
    \vspace{-5mm}
    \caption{The top plot shows the distance between the exact and expected upper bounds for~$\rho$. The bottom plot shows the distance between the corresponding lower bounds.}
    \label{fig:bounds_for_rho}
\end{figure}

We next present simulation results for using the private value of~$\lambda_2$ 
to estimate $d$ and $\rho$. We consider networks of $n=30$ agents with different edge sets 
and hence different values of~$\lambda_2$. We let $\lambda_n=n$  and therefore the 
upper bounds on $d$ and $\rho$ in Theorem~\ref{thm:expectation_bounds_d_rho} 
can reach their worst-case values. We apply the bounded 
Laplace mechanism with $\delta=0.05$ and a range of $\epsilon\in[0.1,2]$.
To illustrate the effects of privacy in bounding diameter, we compute the distance
between the exact (non-private) upper bound on diameter in Lemma~\ref{thm:d_rho_bounds}
and the expected (private) upper bound on diameter in Theorem~\ref{thm:expectation_bounds_d_rho}. 
This distance is shown in the upper plot in Figure~\ref{fig:bounds_for_d}, and the lower
plot shows the analogous distance for the diameter lower bounds. 
Figure~\ref{fig:bounds_for_rho} shows the corresponding upper- and lower-bound distances
for~$\rho$. 

In all plots, we see that the errors induced by privacy are small. Moreover, 
there is a general decrease in the distance between the exact and private bounds as~$\epsilon$ grows. 
Recalling that a larger~$\epsilon$ implies weaker privacy, 
these simulations confirm that weaker privacy guarantees result in smaller differences between the exact and expected bounds for $d$ and $\rho$, as predicted in Remark~\ref{rmk:epsilon_bounds_relation}.

%% file: guidelines_examples.tex
\section{Guidelines and Examples}
\label{sec:examples}
\begin{figure}
    \centering
    \includegraphics[width=.4\textwidth,draft=false]{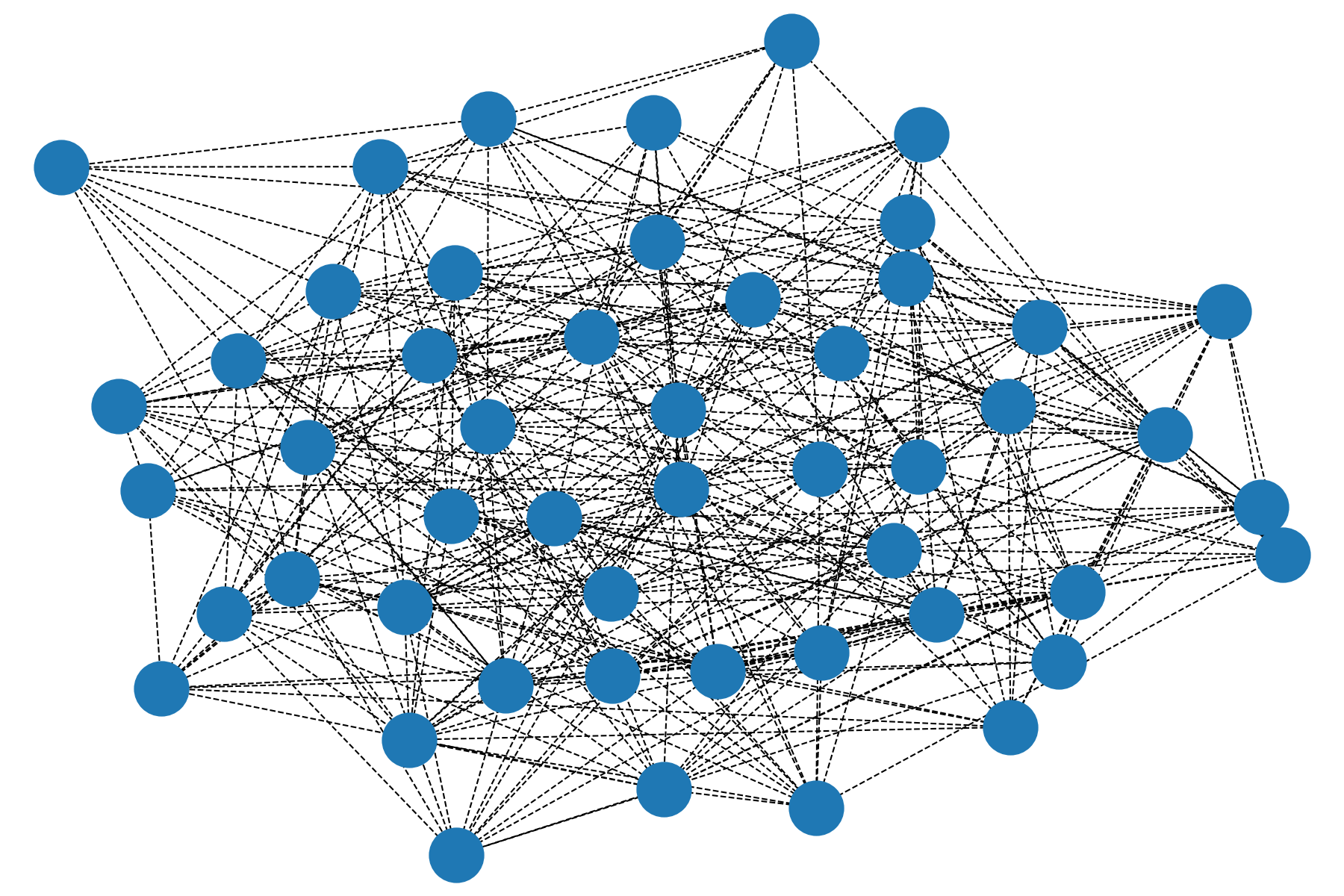}
    \caption{The graph $G$ used in Examples~\ref{ex:example_1}-\ref{ex:kem_example},
    which is an Erdos-Renyi graph with $n=50$ and $p=0.40.$}
    \label{fig:ex1_input_graph}
\end{figure}
In this section, we develop guidelines for providing private responses to queries of the Laplacian eigenvalues, as well as a series of examples to highlight what type of information can be shared via queries of the Laplacian spectrum, 
thereby solving Problem~\ref{prb:prob4}. Recall that a connected graph $G\in\mathcal{G}_n$ has eigenvalues $\lambda_1\leq\lambda_2\leq\dots\leq\lambda_n,$ where $\lambda_1=0$ and $\lambda_2>0.$ In this section, we generate private eigenvalues $\tilde\lambda_i$ according to a mechanism $W_{\lambda_i},$ which we write as $\tilde\lambda_i\sim W_{\lambda_i}.$

The procedure for sharing one private eigenvalue is straightforward. Given a graph $G,$ and privacy parameters $\epsilon$ and $\delta,$ we can compute the eigenvalue $\lambda_i$ and the minimum $b$ required for $(\epsilon,\delta)-$differential privacy, either edge or node, then add noise with the bounded Laplace mechanism to get the private eigenvalue $\tilde\lambda_i$. More care must be taken when answering queries of multiple eigenvalues or the entire spectrum. Specifically, since we only consider connected graphs we will always have $\lambda_1=0$ and thus there is no need to privatize it. Furthermore, for $\lambda_i$ with $i\in\{2,\dots,n\}$ we can define $n-1$ independent mechanisms that provide $(\epsilon,\delta)-$differential privacy to each $\lambda_i.$ In general, since we have $n-1$ queries that are each individually $(\epsilon,\delta)-$differentially private, the privacy level for querying the entire spectrum is $\big((n-1)\epsilon, (n-1)\delta\big)-$differentially private due to the Composition Theorem~\cite[Theorem 3.16]{dwork_algorithmic_2013}. 
After privatizing the spectrum, the set $\{\tilde\lambda_i\}_{i=1}^n$ is no longer guaranteed to have the ordering $\tilde\lambda_1\leq\dots\leq\tilde\lambda_n.$ In applications where the sorting of the private values is critical we can sort the private values prior to sharing them. Sorting does not harm privacy 
because it is post-processing on privatized data, but it will change the statistics of each $\tilde\lambda_i.$

For the remainder of this section we provide a series of examples illustrating the accuracy and utility of the edge privacy mechanism developed in Theorem~\ref{thm:mechanim_variance_ub_edge}. 
In each of the examples, we calculate a metric to quantify accuracy of the private information, and Table~\ref{tab:graphs} gives statistical summaries of these quantities.
\begin{table}
\centering
\begin{tabular}{| c | c | c | c | c |}
    \hline 
     Quantity & $\epsilon$ & $n$ & Average \% error & Variance of Error \\ \hline
    $\lambda_2(G)$  & $0.60$ & $50$ & $8.81\%$ & $0.26$ \\ \hline
    $Tr(L(G))$ & $0.35$ & $50$ & $5.15\%$ & $0.01$ \\ \hline
   $K(P)$ & $1.00$ & $50$ & $4.42\%$ & $0.01$ \\ \hline
    $\phi(G)$ & $2.50$ & $14$ & $9.01\%$ & $0.27$ \\ \hline
  \end{tabular}
  \caption{
Summary of the quantities computed in Section~\ref{sec:examples}. Values were computed using $M=10^4$ private spectrum values.
There are no columns for~$A$ and~$\delta$ since they are fixed at~$A = 2$ and~$\delta = 0.05$ for all simulations. 
  }
  \label{tab:graphs}
 \end{table}

\begin{figure}
    \centering
    \includegraphics[width=.51\textwidth,draft=false]{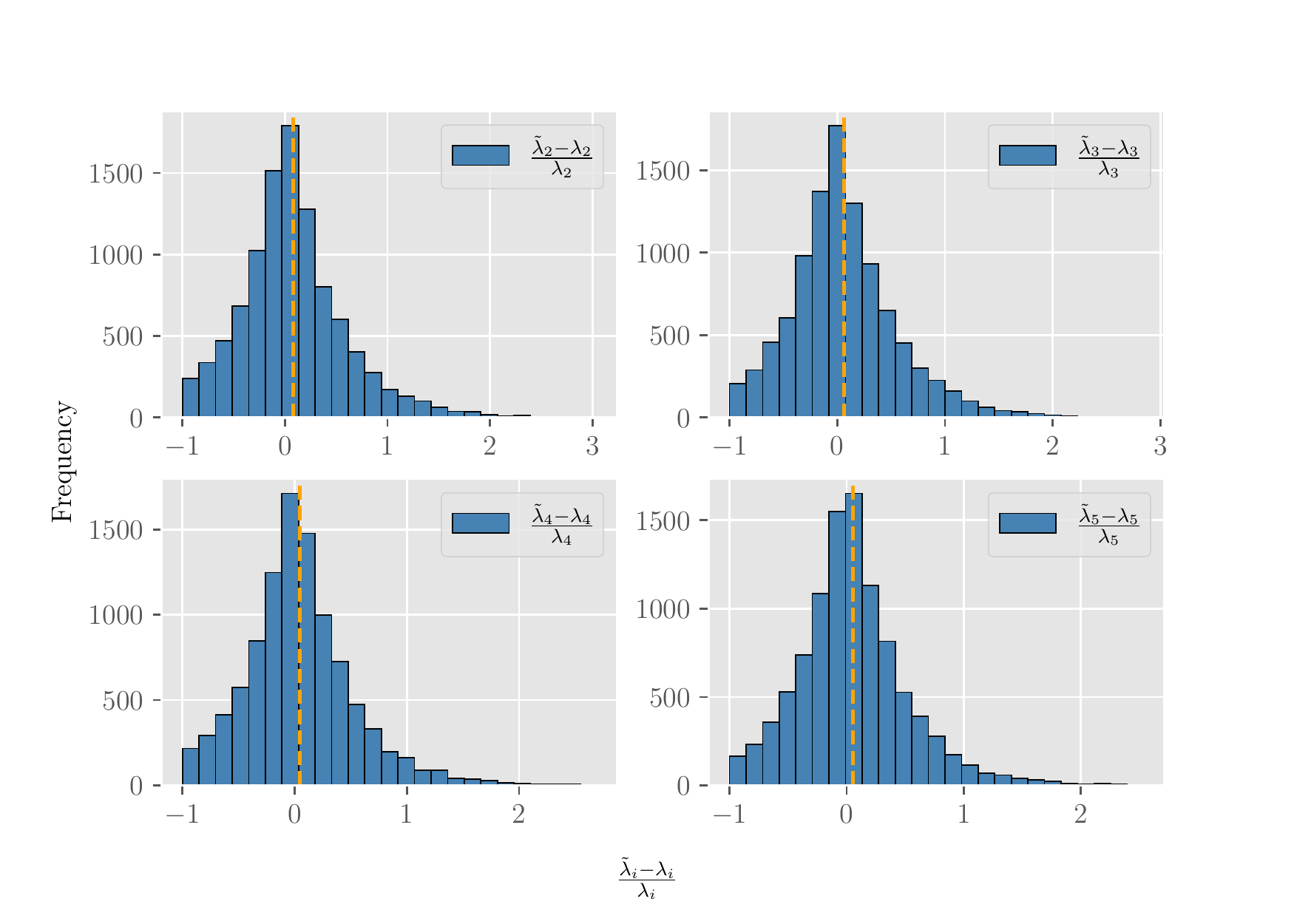}
    \caption{Errors in private values of $\lambda_i$ for $i\in\{2,\dots, 5\}.$ These results illustrate that edge privacy is able to achieve high accuracy,
    even under strong privacy. 
    }
    \label{fig:error_hist}
\end{figure}

\begin{example}[\emph{Accuracy}] \label{ex:example_1}
Fix $G\in\mathcal{G}_{50}$ to be the graph shown in Figure~\ref{fig:ex1_input_graph}. Fix $\epsilon=0.6,\delta=0.05$, and $A=2.$ We generated $M=10^4$ private $\tilde\lambda_i$'s 
for each $i \in \{2,\dots, 5\}$ using an edge privacy mechanism $W_{\lambda_i}^e$ with parameter $b_e$. Solving for the minimum $b_e$ required for $(0.6,0.05)-$differential privacy gives $b_e=6.386.$ To quantify the accuracy of the private spectrum for a fixed $\epsilon$ and $\delta$ we analyze $\tilde\lambda_i-\lambda_i$ for~$i \in \{2, \ldots, 5\}$. A histogram of the accuracy for the $M=10^4$ queries is shown in Figure~\ref{fig:error_hist}. For each of the eigenvalues, the error in the private information is heavily concentrated near $0.$ This trend persists for the rest of the $N=50$ eigenvalues as well as for larger networks with larger values of~$N$. This shows that edge privacy provides accurate spectrum values for large networks, even with strong privacy.

In Figure~\ref{fig:error_hist}, it appears that there is a slight bias in the private spectrum values because the plots are not perfectly symmetric. This bias is made precise by Theorem~\ref{thm:edge_privacy_acc}, and
it is a function of the underlying graph $G$ through its eigenvalues and a function of 
the privacy parameters $\epsilon$ and $\delta$ through $b_e.$ This bias appears as a result of adding bounded noise. Specifically, the density we use to generate $\tilde\lambda_i$ has a peak at the true value $\lambda_i$ but is only supported on the interval $[0,n]$, which means that 
the expected value will not be $\lambda_i$ unless~$\lambda_i = \frac{n}{2}$.  
Nonetheless, 
Figure~\ref{fig:error_hist} shows that this bias is small even when using strong privacy.\hfill$\triangle$

\end{example}
\begin{figure}
    \centering
    \includegraphics[width=.5\textwidth,draft=false]{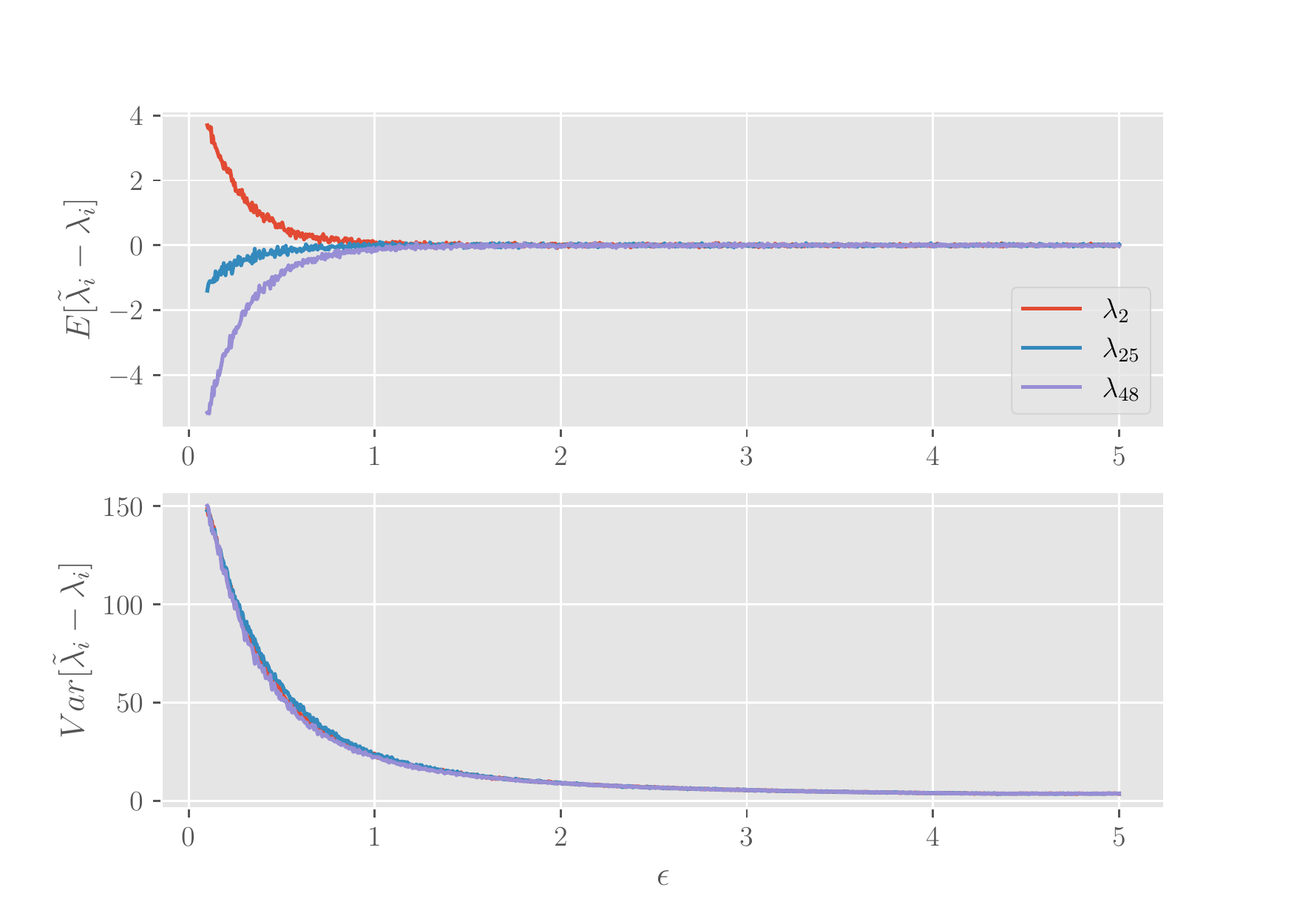}
    \caption{The empirical mean and variance of the error in $\tilde\lambda_i$ for $i\in \{2,25,48\}$ and $\epsilon\in[0.1,5].$ For each $\epsilon,$ $M=10^4$ private values were generated to empirically compute the values of $E[\tilde\lambda_i-\lambda_i]$ and $\textnormal{Var}[\tilde\lambda_i-\lambda_i].$}
    \label{fig:effect_of_eps}
\end{figure}
\begin{example}[\emph{The Effect of $\epsilon$}]
 \label{ex:eps_example}
 
Fix $G\in\mathcal{G}_{50}$ to be the graph shown in Figure~\ref{fig:ex1_input_graph}. Fix $\delta=0.05$ and $A=2.$ Let $\epsilon$ vary and take on values $\epsilon\in[0.1,5].$ Then for each $\epsilon$, generate $M=10^4$ private $\tilde\lambda_i$'s for $i\in \{2, \dots, 5\}$ using an edge privacy mechanism $W_{\lambda_i}^e$ with parameter $b_e$. For a given $\epsilon$ and eigenvalue $\lambda_i,$ we quantify the quality of the private information with the empirical values of $E[\tilde\lambda_i-\lambda_i]$ and $\textnormal{Var}[\tilde\lambda_i-\lambda_i]$ taken over the $M=10^4$ private values. Figure~\ref{fig:effect_of_eps} presents the values of $E[\tilde\lambda_i-\lambda_i]$ and $\textnormal{Var}[\tilde\lambda_i-\lambda_i]$ for $\epsilon\in[0.1,5].$ Recall that a larger $\epsilon$ implies weaker privacy. 

In Figure~\ref{fig:effect_of_eps}, as $\epsilon$ grows and privacy is weakened, both $E[\tilde\lambda_i-\lambda_i]$ and $\textnormal{Var}[\tilde\lambda_i-\lambda_i]$ converge to $0$ relatively quickly. This trend is consistent across the entire spectrum of the graph Laplacian. This shows that even with relatively strong privacy, for example $\epsilon=2$, the private spectra we share are highly accurate. Here we also we see that under strong privacy, given by small $\epsilon,$ we are sharing values of $\lambda_2$ that are much larger than the true value, and we are sharing much smaller values of $\lambda_{48}.$ This occurs because of adding bounded noise and because $\lambda_{2}$ and $\lambda_{48}$ are near the boundaries of the allowable output range $[0,n].$ This example also illustrates the loss of accuracy as privacy is strengthened.
\end{example}

\begin{figure}
    \centering
    \includegraphics[width=.5\textwidth,draft=false]{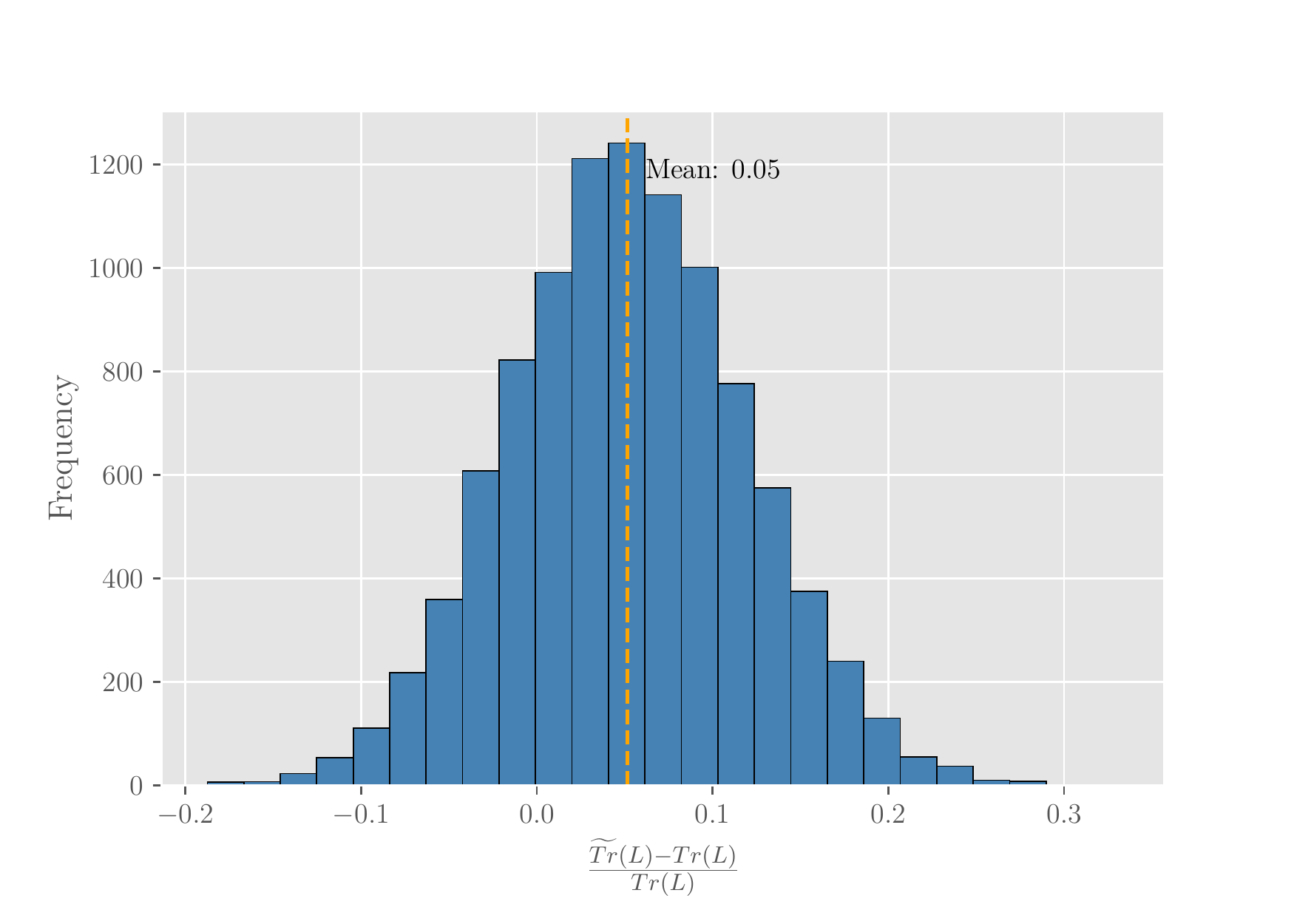}
    \caption{Values of $\widetilde{Tr}(L)$ that are computed using private eigenvalues. For privacy parameters $\epsilon=0.4$ and $\delta=0.05$, $10^4$ sets of eigenvalues were generated.}
    \label{fig:trace_hist}
\end{figure}

\begin{example}[\emph{Trace of the Laplacian}]
\label{ex:trace}
Fix $G\in\mathcal{G}_{50}$ to be the graph shown in Figure~\ref{fig:ex1_input_graph}. Fix $\epsilon=0.4$, $\delta=0.05$, and $A=2$. Recall that the trace of a matrix $R\in\mathbb{R}^{n\times n}$ is given by the sum of its eigenvalues, i.e., $Tr(R)=\sum_{i=1}^n \lambda_i(R).$ Applying this to the graph Laplacian, we have $Tr(L)=\sum_{i=1}^n \lambda_i(L).$ The trace of the graph Laplacian can, for example, be used to compute the average degree of the network as $d_{avg} = \frac{Tr(L)}{n}.$ Suppose that we do not have access to $G$ or $L(G)$ and we only have the private spectrum values $\{\tilde\lambda_i\}_{i=1}^N$. Then we can use these eigenvalues to estimate the trace of $L$ as $\widetilde{Tr}(L)=\sum_{i=1}^n \tilde\lambda_i(L).$ To analyze the accuracy of this estimate, $M=10^{4}$ sets of private spectra were generated and used to estimate the trace. In Figure~\ref{fig:trace_hist} we give a histogram of values of  $\widetilde{Tr}(L)-Tr(L)$ for these trace estimates. The trace of the graph appearing in Figure~\ref{fig:ex1_input_graph} is $Tr(L)=736$ and the average estimate over the $M=10^4$ queries was $774.$ We can see in Figure~\ref{fig:trace_hist} that edge privacy generally provides accurate estimates of the trace, with
the majority of private trace estimates falling within~$\pm10\%$ of the
true trace value. 

However, there is a bias in the distribution of private trace estimates, 
and we tend to overestimate the trace. To quantify this overestimate, we analyze $E\left[\widetilde{Tr}(L)-Tr(L)\right].$ 
Plugging in $Tr(L)=\sum_{i=1}^n \lambda_i(L)$ and simplifying gives
\begin{equation}
E\left[\widetilde{Tr}(L)-Tr(L)\right]=\sum_{i=1}^{n}E\left[\tilde{\lambda}_{i}(L)\right]-\lambda_{i}(L).
\end{equation}
Then applying Theorem~\ref{thm:edge_privacy_acc} gives
\begin{equation}
    E\left[\widetilde{Tr}(L)-Tr(L)\right]=\sum_{i=1}^{n}\frac{1}{2C(\lambda_{i},b)}\left(2\lambda_{i}+be^{-\frac{\lambda_{i}}{b}} \right. \left.-(n+b)e^{-\frac{n-\lambda_{i}}{b}}\right)-\lambda_{i}(L),
\end{equation}
where~$C$ is from Definition~\ref{dfn:bounded_laplace}. 
In Example 1, there was a small bias in the values of $\tilde\lambda_i$ due to using bounded noise to achieve differential privacy. Here the bias for the trace is larger 
because we are summing each $\tilde\lambda_i$ and the bias is amplified due to summing
biased terms.
Nonetheless, accurate trace estimates can still be attained, even under strong privacy. 
\hfill$\triangle$
\end{example}

\begin{figure}
    \centering
    \includegraphics[width=.5\textwidth,draft=false]{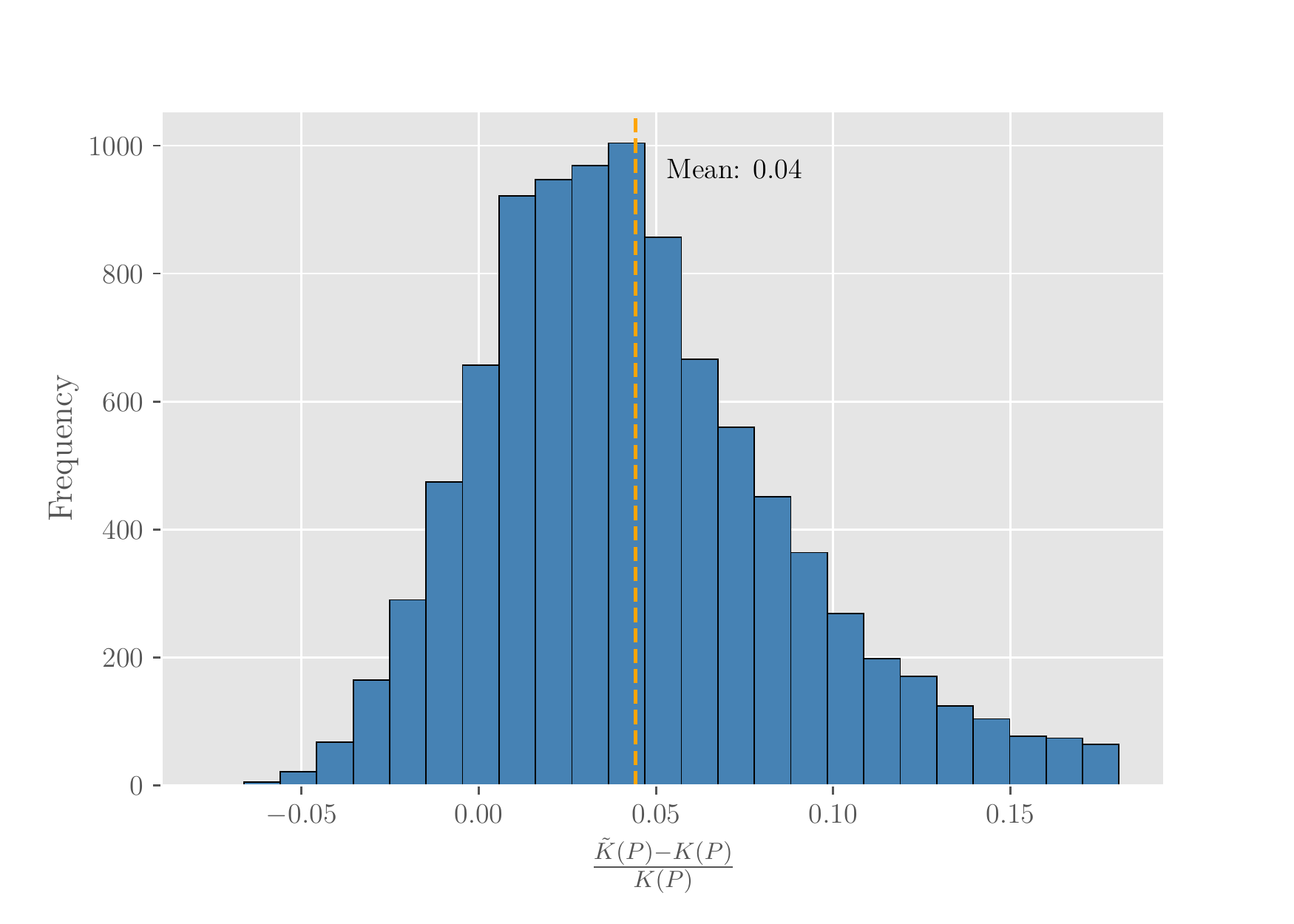}
    \caption{The values of $\frac{\widetilde{K}(P)-K(P)}{K(P)}$ in Example~\ref{ex:kem_example}. 
    The average value of $\frac{\widetilde{K}(P)-K(P)}{K(P)}$ is $0.2519$.
    In general, we overestimate the Kemeny constant of the graph $G$ in Figure~\ref{fig:ex1_input_graph}, but the majority
    of results are within~$\pm10\%$ of its true value.}
    \label{fig:kem_hist}
\end{figure}
\begin{example}[\emph{Kemeny's Constant}]\label{ex:kem_example}
In network control, network level discrete-time consensus dynamics are governed by the matrix $P=I-\gamma L(G),$ where $\gamma$ is a step-size which must obey $\gamma\leq\frac{1}{\max_i d_i}$ in order to achieve 
consensus\cite[Theorem 2]{olfati2007consensus}. 
When $G$ is a connected, undirected graph, $P$ can be interpreted as the transition matrix of a symmetric Markov chain. 
The Kemeny constant of a Markov chain is the expected time it takes to transition from a state $i$ in a Markov chain to another state sampled from its stationary distribution and can be used to compute the error in consensus protocols subject to noise \cite{jadbabaie2018scaling}.
The Kemeny constant of the Markov chain with transition matrix 
$P=I-\gamma L(G)$ can be computed as $K(P)=\sum_{i=2}^n \frac{1}{1-\lambda_i(P)}$ \cite{levene2002kemeny}. Note that $\lambda_i(P)=1-\gamma\lambda_i(L)$ and thus $K(P)=\frac{1}{\gamma}\sum_{i=2}^n \frac{1}{\lambda_i(L)}.$ Given private spectrum values we can estimate the Kemeny constant as $\widetilde{K}(P)=\frac{1}{\gamma}\sum_{i=2}^n \frac{1}{\tilde\lambda_i(L)}.$ We fix $\gamma=\frac{1}{N}$, and
with this step-size the graph in Figure~\ref{fig:ex1_input_graph} has $K(P)=102.70.$ 

We now fix $\epsilon=1.0$, $\delta=0.05$, and $A=2$. We generate $M=10^4$ private spectra for $G$ in Figure~\ref{fig:ex1_input_graph} and these values are used to compute $\widetilde{K}(P).$
To quantify the accuracy of the estimates of the Kemeny constant, we analyze the relative error
$\frac{\widetilde{K}(P)-K(P)}{K(P)}$ whose values for the $M=10^4$ private spectra are presented in Figure~\ref{fig:kem_hist}. Here, we can see that we overestimate the Kemeny constant,
but the average error for these queries is only $4.42\%.$ This shows that sharing the private spectrum can share relatively accurate information about the Kemeny constant and thus about discrete-time consensus dynamics while providing edge differential privacy. 

\end{example}

\begin{figure}
    \centering
    \includegraphics[width=.3\textwidth,draft=false]{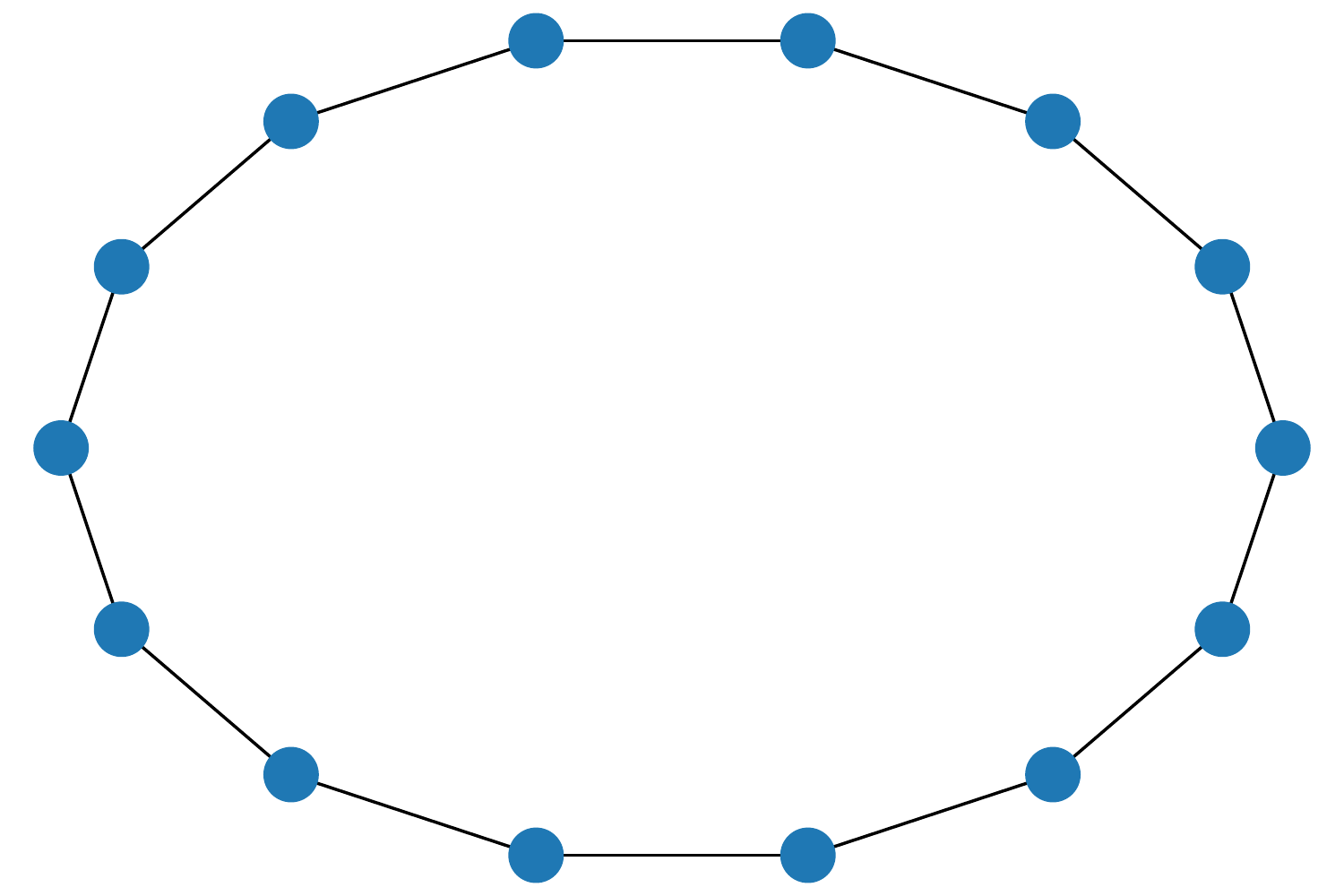}
    \caption{The cycle graph on $n=14$ nodes, $C_{14},$ used in Example~\ref{ex:cheeger}.}
    \label{fig:cycle}
\end{figure}
\begin{figure}
    \centering
    \includegraphics[width=.5\textwidth,draft=false]{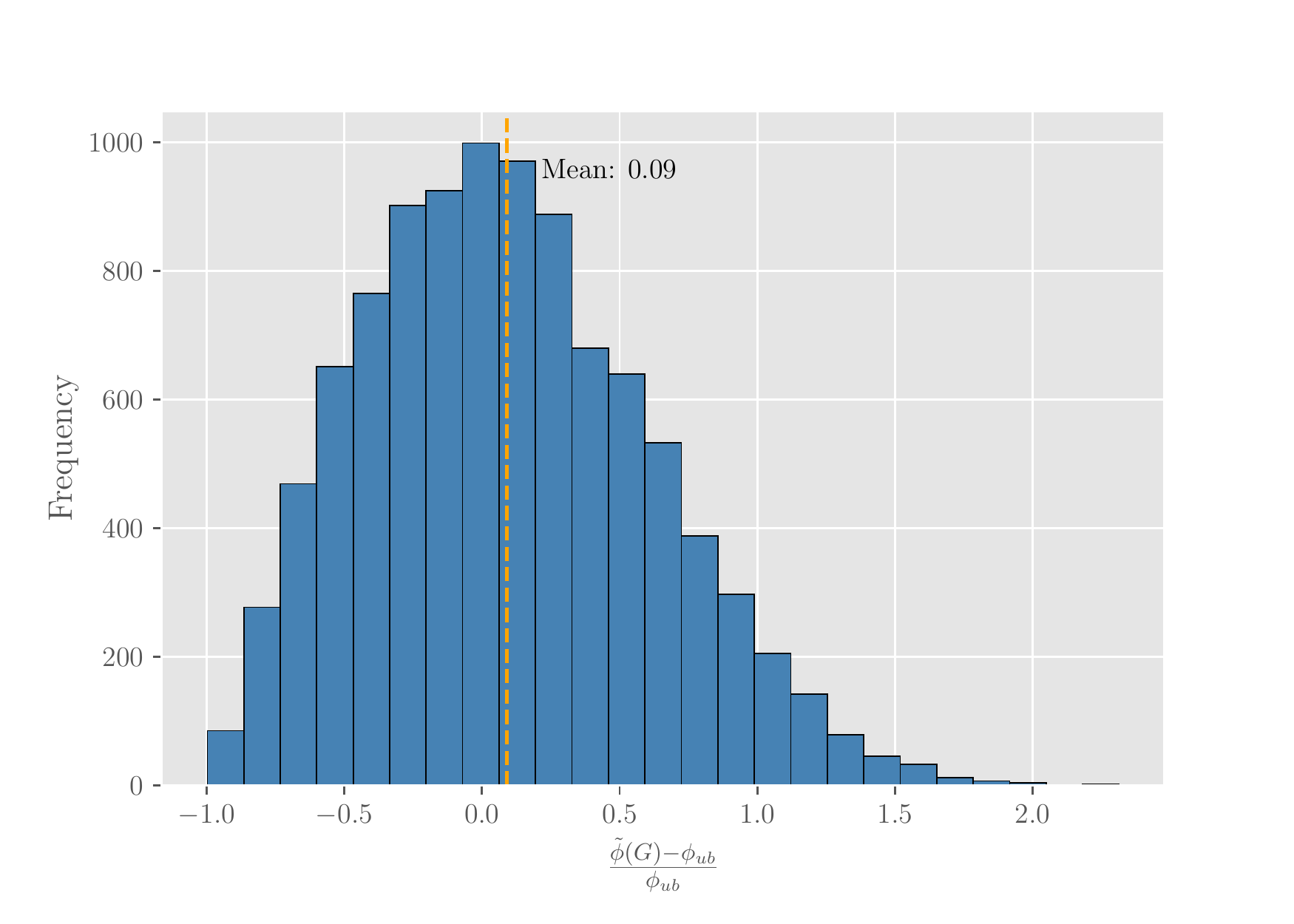}
    \caption{The error in the estimate of Cheeger's constant for the cycle graph in Figure~\ref{fig:cycle} and the parameters in Example~\ref{ex:cheeger}. We usually overestimate Cheeger's constant using private information. This means we are estimating that the graph is more connected than it really is, though these estimates are often fairly accurate. 
    }
    \label{fig:cheeger_hist}
\end{figure}
\begin{example}[\emph{Cheeger's Inequality}]
\label{ex:cheeger}
In this example we discuss how private Laplacian spectra 
can be used to estimate the isoperimetric number, $\phi(G),$ of a graph $G$. 
The isoperimetric number, or the Cheeger constant, is a measure of how connected a graph is or more specifically how easy it is to disconnect a graph~\cite{mohar1989isoperimetric}.
In general, the isoperimetric number is NP-hard to compute and Cheeger's inequality gives an easily computable upper bound on the isoperimetric number via
$\phi(G)\leq\phi_{ub} := \sqrt{\lambda_2\left(\max_i d_i-\lambda_2\right)}$ \cite[Theorem 4.2]{mohar1989isoperimetric}.

In this example, we estimate $\phi(G)$ using Cheeger's inequality
for cases in which
we do not have access to $G$ and only have its private Laplacian spectrum. 
To estimate $\lambda_2,$ we use the private value $\tilde\lambda_2.$ For $\max_i d_i$, we estimate this with $\tilde{d}(G)=\frac{1}{n}\sum_{i=1}^n\tilde\lambda_i.$ Then plugging these estimates into Cheeger's inequality, we have the estimate
$\tilde\phi(G)=\sqrt{\tilde\lambda_2\left(2\tilde d(G)-\tilde\lambda_2\right)}.$

Since the isoperimetric number is not feasible to compute for large networks, we cannot run simulations on the graph appearing in Figure~\ref{fig:ex1_input_graph} to demonstrate the accuracy of our estimates. Thus, we fix $G$ to be the cycle or ring graph on $n$ nodes, $C_n,$ which has 
a known Cheeger's constant of $\phi(C_n)=\frac{4}{n}$ \cite{godsil2001algebraic}. For this example, we fix $n=14$. 
The graph $G = C_{14}$ is shown in Figure~\ref{fig:cycle} and $\phi(C_{14})=\frac{2}{7}$. To analyze the accuracy of 
using Cheeger's inequality with private spectra, 
we generate $M=10^4$ private Laplacian spectra $\{\tilde\lambda_i\}_{i=1}^n$ and use them to privately estimate~$\phi(G)$. 

Before discussing the accuracy of our estimates we will discuss the accuracy of the Cheeger's inequality itself. For $C_{14},$ we have $\phi(C_{14})=\frac{2}{7}=0.2857$ and plugging in $\lambda_2$ and $\max_id_i = 2$ into Cheeger's inequality gives $\phi_{ub}=0.8678.$ This is more than $3$ times the true value. Thus,
to distinguish between errors inherent to Cheeger's inequality itself and errors due to privacy, 
we will compare our estimate to the upper bound from Cheeger's inequality, $\phi_{ub}$. 

In Figure~\ref{fig:cheeger_hist}, we show the accuracy of the resulting estimates given by $\frac{\tilde\phi(G)-\phi_{ub}}{\phi_{ub}}$ for $10^4$ queries satisfying $(2.5,0.05)-$differential privacy. Here, we typically over estimate Cheeger's constant. This means that we are estimating that the graph is more connected than it truly is. Comparing to the non-private Cheeger's inequality upper bound given by~$\phi_{ub}$, the use of private spectra in computations results in a slightly looser bound on average. However the estimates are relatively accurate with an average normalized error 
of $9.01\%$, with a variance of only $0.27.$ 
Overall, this example shows that using private spectrum information to estimate the isoperimetric number is relatively accurate and does not have much more error than when true spectrum values are used.
\hfill$\triangle$


\end{example}

%% file: appendix.tex
\appendix
\subsection{Proof of Lemma~1} \label{apdx:edge_sens}
Fix an adjacency parameter~$A \in \N$ and
consider two graphs $G, G' \in \mathcal{G}_n$ such that $\adj_{e,A}(G, G')=1$. Denote their corresponding graph Laplacians by $L$ and $L'$, and
define the matrix $P$ such that $L' = L+P$. Then, we write 
\begin{align}
 \Delta\lambda_{e,i} &=  \max_{G,G' \in \mathcal{G}_n} \left|\lambda_{i}\left(L+P\right)-\lambda_{i}\left(L\right)\right|. 
\end{align}
We will use the following lemma. 
\begin{lemma}[$\!\!\!${\cite[Theorem 8.4.11]{bernstein2009matrix}}]\label{lem:eig_perturbation} Let $A,B\in\mathbb{R}^{n\times n}$ be two symmetric matrices. Then $\lambda_i(A+B)\leq\lambda_{i+j}(A)+\lambda_{n-j}(B).$
\end{lemma}

Applying Lemma~\ref{lem:eig_perturbation} to split up~$\lambda_i(L+P$), we obtain $\Delta\lambda_{i,e} \le\lambda_{i}(L)+\lambda_{n}(P) - \lambda_{i}(L) = \lambda_n(P). $
The matrix~$P$ encodes the differences between~$L$ and~$L'$ as follows. For any~$i$, if the diagonal
entry~$P_{ii} = 1$, then node~$i$ has one more edge in~$G'$ than it does in~$G$. If~$P_{ii} = -1$, then
node~$i$ has one fewer edge in~$G'$ than it does in~$G$.
Other values of~$P_{ii}$ indicate the addition or removal of more edges. 
Because edge adjacency allows for the addition or removal of up to~$A$
edges, we have~$-A \leq P_{ii} \leq A$. 

For off-diagonal entries,~$P_{ij} = 1$ 
indicates that~$G'$ contains the edge~$(i, j)$ and~$G$ does not; the converse holds
if~${P_{ij} = -1}$. Then, for any row of~$P$, the diagonal entry has absolute value at most~$A$,
and the absolute sum of the off-diagonal entries is at most~$A$. 
By Ger\v{s}gorin's circle theorem \cite[Fact 4.10.16.]{bernstein2009matrix}, we have $\lambda_n(P)\le 2A.$\hfill$\blacksquare$

\subsection{Proof of Lemma~2}
\label{apdx:node_sens}
Let $G^{-}\in\mathcal{G}^{n-1}$
be the graph obtained by deleting vertex $v$ and its incident edges
from $G$. Then
\begin{equation}
-1\leq\lambda_{2}(G^{-})-\lambda_{2}(G)\leq n-2.\label{node_del_bound}
\end{equation}
The lower bound in Equation~\eqref{node_del_bound} is
given in \cite[Theorem 1.1]{kirkland2010algebraic} and the upper bound in \cite[Theorem 2.3]{kirkland2010algebraic}. 

We now analyze the case where we add a node $v$ and arbitrary incident edges to obtain the graph $G^+\in\mathcal{G}^{n+1}.$ Using \cite[Theorem 1.1]{kirkland2010algebraic}, we have that 
\begin{equation} \label{eq:node_sens_1}
\lambda_{2}(G^{+})-\lambda_{2}(G)\leq 1.
\end{equation}
To lower bound $\lambda_{2}(G^{+})-\lambda_{2}(G),$ we can apply the same methods to obtain the upper bound in Equation~\eqref{node_del_bound}. Let $K_n$ be the complete graph on $n$ nodes. Suppose that $G\neq K_{n}.$ Then $\lambda_{2}(G)\leq n-2$ and
\begin{equation}
\lambda_{2}(G^+)-\lambda_{2}(G)\geq\lambda_{2}(G^+)-(n+2),
\end{equation}
and since $\lambda_{2}(G^+)>0$ we have
\begin{equation} \label{eq:node_sens_2}
\lambda_{2}(G^+)-\lambda_{2}(G) > 2-n
\end{equation}
for~$G \neq K_n$. 

Now suppose $G=K_{n}$ and we add a node with degree $d\geq1,$
i.e., we add $d$ incident edges. It can be shown that $\lambda_{2}(G^+)=d$ \cite[Theorem 2.3]{kirkland2010algebraic},
and thus
\begin{align}
\lambda_{2}(G^+)-\lambda_{2}(G) &=    d-n\\
                                &\geq 1-n. \label{eq:node_sens_3}
\end{align}
With this, we achieve equality when $d=1,$ which occurs when the added node has one incident edge. Thus, for any $G$ and $G^+$ obtained by adding a node,
we combine~\eqref{eq:node_sens_1}, \eqref{eq:node_sens_2}, and~\eqref{eq:node_sens_3} to find
\begin{equation}
1-n\leq\lambda_{2}(G^+)-\lambda_{2}(G)\leq 1.\label{node_add_bound}
\end{equation}
Lastly, Equations~\eqref{node_del_bound} and~\eqref{node_add_bound} imply that for any $G, G'$ satisfying $\textnormal{Adj}_{n}(G,G')=1$, 
we have
\begin{equation}
|\lambda_{2}(G')-\lambda_{2}(G)|\leq\max\left\{ n-1,n-2,1\right\},
\end{equation}
and thus $\Delta\lambda_{2,n}\leq n-1.$\hfill$\blacksquare$

\subsection{Proof of Corollary~1} \label{appdx:node_weak}
We begin with the condition in Theorem~2 that $$b_n\geq\frac{n-1}{\epsilon-\log\left(\frac{2-e^{-\frac{n-1}{b_n}}-e^{-\frac{1}{b_n}}}{1-e^{-\frac{n}{b_n}}}\right)-\log(1-\delta)}.$$ To find a necessary condition, we must find a lower bound on the right hand side of the expression. To achieve this, we focus on the argument of the log first. We note that if 
\begin{equation}
\frac{2-e^{-\frac{n-1}{b_n}}-e^{-\frac{1}{b_n}}}{1-e^{-\frac{n}{b_n}}}>1,
\end{equation}
then the $\log$ term will be positive and we can eliminate it from the denominator to obtain the desired bound. 

To show this positivity we begin with the fact that $1>e^{-\frac{1}{b_n}}$ and $1-e^{-\frac{n-1}{b_n}}>0$ for any $n\geq 2$ and $b_n>0$.  Dividing the first inequality by $1-e^{-\frac{n-1}{b_n}}$ gives
\begin{align}
    \frac{1}{1-e^{-\frac{n-1}{b_n}}}&>\frac{e^{-\frac{1}{b_n}}}{1-e^{-\frac{n-1}{b_n}}}.\\
\end{align}
We rearrange this and manipulate the inequality further to find 
\begin{align}
    1-e^{-\frac{n-1}{b_n}}	&> e^{-\frac{1}{b_n}}(1-e^{-\frac{n-1}{b_n}})\\
    1-e^{-\frac{n-1}{b_n}}	&> e^{-\frac{1}{b_n}}-e^{-\frac{n}{b_n}}\\
    1-e^{-\frac{n-1}{b_n}}-e^{-\frac{1}{b_n}}	&>-e^{-\frac{n}{b_n}}\\
    2-e^{-\frac{n-1}{b_n}}-e^{-\frac{1}{b_n}}	&>1-e^{-\frac{n}{b_n}}\\
    \frac{2-e^{-\frac{n-1}{b_n}}-e^{-\frac{1}{b_n}}}{1-e^{-\frac{n}{b_n}}}	&>1.
\end{align}
This implies that $-\log\left(\frac{2-e^{-\frac{n-1}{b_n}}-e^{-\frac{1}{b_n}}}{1-e^{-\frac{n}{b_n}}}\right)<0$, which gives us $$\frac{n-1}{\epsilon-\log\left(\frac{2-e^{-\frac{n-1}{b_n}}-e^{-\frac{1}{b_n}}}{1-e^{-\frac{n}{b_n}}}\right)-\log(1-\delta)} > \frac{n-1}{\epsilon-\log(1-\delta)}.$$ Thus, choosing $b_n$ according to $b_n>\frac{n-1}{\epsilon-\log(1-\delta)}$ is necessary to satisfy $(\epsilon,\ \delta)$-node differential privacy.
\hfill$\blacksquare$

\subsection{Proof of Corollary~2}
\label{appdx:edge_weak}
Following the proof of Corollary~1, we show that the argument of the $\log$ 
in Theorem~1 is larger than $1.$ We begin with $1> e^{-\frac{2A}{b}}$ and follow a similar 
sequence of steps to~Corollary~1: 
\begin{align}
    1&> e^{-\frac{2A}{b}}\\
    1-e^{-\frac{n-2A}{b}}&> e^{-\frac{2A}{b}}\left(1-e^{-\frac{n-2A}{b}}\right)\\
    1-e^{-\frac{n-2A}{b}}&> e^{-\frac{2A}{b}}-e^{-\frac{n}{b}}\\
    2-e^{-\frac{2A}{b}}-e^{-\frac{n-2A}{b}}&>1-e^{-\frac{n}{b}}\\
    \frac{2-e^{-\frac{2A}{b}}-e^{-\frac{n-2A}{b}}}{1-e^{-\frac{n}{b}}}&>1.
\end{align}
Thus, by an argument similar to that in Corollary~1, we find that 
\begin{equation} 
    b_e>\frac{2A}{\epsilon-\log(1-\delta)}
\end{equation}
is a necessary condition for the satisfaction of~$(\epsilon, \delta)$-edge differential privacy.
\hfill$\blacksquare$

\subsection{Proof of Theorem~3}
\label{apdx:proof_of_acc}
 We first compute $E[\tilde\lambda_i]$ as
\begin{align}
    E\left[\tilde{\lambda}_{i}\right]&=\frac{1}{C(\lambda_{i},b_e)}\frac{1}{2b_e}\int_{0}^{n}xe^{-\frac{\left|x-\lambda_{i}\right|}{b_e}}dx\\&=\frac{1}{C(\lambda_{i},b_e)}\frac{1}{2b_e}\left(\int_{0}^{\lambda_{i}}xe^{-\frac{\lambda_{i}-x}{b_e}}dx+\int_{\lambda_{i}}^{n}xe^{-\frac{x-\lambda_{i}}{b_e}}dx\right)\\&=\frac{1}{2C(\lambda_{i},b_e)}\left(2\lambda_{i}+b_e e^{-\frac{\lambda_{i}}{b_e}}-b_e e^{-\frac{n-\lambda_{i}}{b_e}}-ne^{-\frac{n-\lambda_{i}}{b_e}}\right)\\
    &=\frac{1}{2C(\lambda_{i},b_e)}\left(2\lambda_{i}+b_e e^{-\frac{\lambda_{i}}{b_e}}-(n+b_e)e^{-\frac{n-\lambda_{i}}{b_e}}\right).\label{expected_lambda_i}
\end{align}
With~\eqref{expected_lambda_i} we can compute $E\left[\tilde{\lambda}_{i}-\lambda_{i}\right]$ as
\begin{align}
  E\left[\tilde{\lambda}_{i}-\lambda_{i}\right]	&=E\left[\tilde{\lambda}_{i}\right]-\lambda_{i}
	\\&=\frac{1}{2C(\lambda_{i},b_e)}\left(2\lambda_{i}+b_e e^{-\frac{\lambda_{i}}{b_e}}-(n+b_e)e^{-\frac{n-\lambda_{i}}{b_e}}\right)-\lambda_{i}.
\end{align}\hfill$\blacksquare$
\subsection{Proof of Theorem~4}
\label{sec:proof_bo_bound}
Since both lower bounds are convex functions with respect to $\lambda_2>0$, we can use Jensen's inequality and we have
    \begin{align}
        &E[\tilde{d}] \geq E[\underline{d}(\tilde{\lambda}_2)] = E\left[\frac{4}{n\tilde{\lambda}_2}\right] \geq \frac{4}{nE[\tilde{\lambda}_2]}\\
        &\begin{multlined}E[\tilde{\rho}] \geq E[\underline{\rho}(\tilde{\lambda_2})] = E\left[\frac{2}{(n-1)\tilde{\lambda}_2} + \frac{n-2}{2(n-1)}\right] \\
        \geq \frac{2}{(n-1)E[\tilde{\lambda}_2]} + \frac{n-2}{2(n-1)}.
        \end{multlined}
    \end{align}
    The value of $E[\tilde{\lambda}_2]$ can be computed as
    \begin{align*}
        E[\tilde{\lambda}_2] &= \frac{1}{C_{\lambda_2}(b)}\frac{1}{2b}\int_0^n x e^{-\frac{|x-\lambda_2|}{b}}dx \\
        &=\frac{1}{C_{\lambda_2}(b)}\frac{1}{2b}\left(\int_0^{\lambda_2}xe^{-\frac{\lambda_2-x}{b}}dx + \int_{\lambda_2}^n xe^{-\frac{x-\lambda_2}{b}}dx \right)\\
        &= \frac{1}{2C_{\lambda_2}(b)}\left(2\lambda_2+be^{-\frac{\lambda_2}{b}}-be^{-\frac{n-\lambda_2}{b}}-ne^{-\frac{n-\lambda_2}{b}}\right).
    \end{align*}
    We next compute the expectation term $E\left[\frac{1}{\sqrt{\tilde{\lambda}_2}}\right]$ as
    
    \begin{align*}
        E\left[\frac{1}{\sqrt{\tilde{\lambda}_2}}\right] &= \frac{1}{C_{\lambda_2}(b)}\frac{1}{2b}\int_0^n \frac{1}{\sqrt{x}} e^{-\frac{|x-\lambda_2|}{b}}dx \\
        &=\frac{1}{C_{\lambda_2}(b)}\frac{1}{2b}\left(\int_0^{\lambda_2}\frac{1}{\sqrt{x}} e^{-\frac{\lambda_2-x}{b}}dx + \int_{\lambda_2}^n\frac{1}{\sqrt{x}} e^{-\frac{x-\lambda_2}{b}}dx \right)\\
        &= \frac{1}{C_{\lambda_2}(b)}\frac{1}{2b} \left(\sqrt{\pi}\sqrt{b}e^{-\frac{\lambda_2}{b}}\left(\textrm{erfi}\left(\sqrt{\frac{\lambda_2}{b}}\right)\right)+\sqrt{b}e^{\frac{\lambda_2}{b}}\left(\Gamma\left(\frac{1}{2},\frac{n}{b}\right)-\Gamma\left(\frac{1}{2},\frac{\lambda_2}{b}\right)\right)\right).
    \end{align*}
    Then we can find the desired upper bounds by applying the linearity of expectation.\hfill $\blacksquare$

%% file: main.bbl
\begin{thebibliography}{10}
\providecommand{\url}[1]{#1}
\csname url@samestyle\endcsname
\providecommand{\newblock}{\relax}
\providecommand{\bibinfo}[2]{#2}
\providecommand{\BIBentrySTDinterwordspacing}{\spaceskip=0pt\relax}
\providecommand{\BIBentryALTinterwordstretchfactor}{4}
\providecommand{\BIBentryALTinterwordspacing}{\spaceskip=\fontdimen2\font plus
\BIBentryALTinterwordstretchfactor\fontdimen3\font minus
  \fontdimen4\font\relax}
\providecommand{\BIBforeignlanguage}[2]{{%
\expandafter\ifx\csname l@#1\endcsname\relax
\typeout{** WARNING: IEEEtran.bst: No hyphenation pattern has been}%
\typeout{** loaded for the language `#1'. Using the pattern for}%
\typeout{** the default language instead.}%
\else
\language=\csname l@#1\endcsname
\fi
#2}}
\providecommand{\BIBdecl}{\relax}
\BIBdecl

\bibitem{Ren2005}
{Wei Ren}, R.~W. {Beard}, and E.~M. {Atkins}, ``A survey of consensus problems
  in multi-agent coordination,'' in \emph{Proceedings of the 2005, American
  Control Conference, 2005.}, 2005.

\bibitem{Scott1988}
J.~Scott, ``Social network analysis,'' \emph{Sociology}, vol.~22, no.~1, pp.
  109--127, 1988.

\bibitem{SHIRLEY2005287}
M.~D. Shirley and S.~P. Rushton, ``The impacts of network topology on disease
  spread,'' \emph{Eco. Complexity}, vol.~2, no.~3, pp. 287--299, 2005.

\bibitem{zheng_consensus_2011}
Y.~Zheng, L.~Wang, and Y.~Zhu, ``Consensus of heterogeneous multi-agent
  systems,'' vol.~5, no.~16, pp. 1881--1888.

\bibitem{Mieghem2009}
P.~{Van Mieghem}, J.~{Omic}, and R.~{Kooij}, ``Virus spread in networks,''
  \emph{IEEE/ACM Transactions on Networking}, vol.~17, no.~1, pp. 1--14, 2009.

\bibitem{freitas2020evaluating}
S.~Freitas and D.~H. Chau, ``Evaluating graph vulnerability and robustness
  using tiger,'' 2020.

\bibitem{kasiviswanathan13}
S.~P. Kasiviswanathan, K.~Nissim, S.~Raskhodnikova, and A.~Smith, ``Analyzing
  graphs with node differential privacy,'' in \emph{Proceedings of the 10th
  Theory of Cryptography Conference on Theory of Cryptography}.\hskip 1em plus
  0.5em minus 0.4em\relax Springer-Verlag, 2013, p. 457–476.

\bibitem{karwa14}
V.~Karwa, S.~Raskhodnikova, A.~Smith, and G.~Yaroslavtsev, ``Private analysis
  of graph structure,'' \emph{ACM Trans. Database Syst.}, vol.~39, no.~3, 2014.

\bibitem{kasiviswanathan14}
S.~P. Kasiviswanathan and A.~Smith, ``On the ’semantics’ of differential
  privacy: A bayesian formulation,'' \emph{Journal of Privacy and
  Confidentiality}, vol.~6, no.~1, Jun. 2014.

\bibitem{dwork_algorithmic_2013}
C.~Dwork and A.~Roth, ``The algorithmic foundations of differential privacy,''
  vol.~9, no.~3, pp. 211--407.

\bibitem{day16}
W.-Y. Day, N.~Li, and M.~Lyu, ``Publishing graph degree distribution with node
  differential privacy,'' in \emph{Proceedings of the 2016 International
  Conference on Management of Data}, 2016, p. 123–138.

\bibitem{shen13}
E.~Shen and T.~Yu, ``Mining frequent graph patterns with differential
  privacy,'' in \emph{Proceedings of the 19th ACM International Conference on
  Knowledge Discovery and Data Mining}, 2013, pp. 545--553.

\bibitem{Ding2018}
X.~Ding, X.~Zhang, Z.~Bao, and H.~Jin, ``Privacy-preserving triangle counting
  in large graphs,'' in \emph{Proceedings of the 27th ACM International
  Conference on Information and Knowledge Management}.\hskip 1em plus 0.5em
  minus 0.4em\relax Association for Computing Machinery, 2018, p. 1283–1292.

\bibitem{Hay2009}
M.~{Hay}, C.~{Li}, G.~{Miklau}, and D.~{Jensen}, ``Accurate estimation of the
  degree distribution of private networks,'' in \emph{2009 Ninth IEEE
  International Conference on Data Mining}, 2009, pp. 169--178.

\bibitem{Task2012}
C.~{Task} and C.~{Clifton}, ``A guide to differential privacy theory in social
  network analysis,'' in \emph{International Conference on Advances in Social
  Networks Analysis and Mining}, 2012, pp. 411--417.

\bibitem{abreu07}
N.~M.~M. {de Abreu}, ``Old and new results on algebraic connectivity of
  graphs,'' \emph{Linear Algebra and its Applications}, vol. 423, no.~1, pp.
  53--73, 2007.

\bibitem{karwa2011private}
V.~Karwa, S.~Raskhodnikova, A.~Smith, and G.~Yaroslavtsev, ``Private analysis
  of graph structure,'' \emph{Proceedings of the VLDB Endowment}, vol.~4,
  no.~11, pp. 1146--1157, 2011.

\bibitem{hawkins2020differentially}
C.~Hawkins and M.~Hale, ``Differentially private formation control,'' in
  \emph{2020 59th IEEE Conference on Decision and Control (CDC)}, 2020.

\bibitem{gohari2020privacy}
P.~Gohari, M.~Hale, and U.~Topcu, ``Privacy-preserving policy synthesis in
  markov decision processes,'' in \emph{2020 59th IEEE Conference on Decision
  and Control (CDC)}, 2020.

\bibitem{gohari2021differential}
P.~Gohari, B.~Wu, C.~Hawkins, M.~Hale, and U.~Topcu, ``Differential privacy on
  the unit simplex via the dirichlet mechanism,'' \emph{IEEE Transactions on
  Information Forensics and Security}, vol.~16, 2021.

\bibitem{cortes2016differential}
J.~Cort{\'e}s, G.~E. Dullerud, S.~Han, J.~Le~Ny, S.~Mitra, and G.~J. Pappas,
  ``Differential privacy in control and network systems,'' in \emph{2016 IEEE
  55th Conference on Decision and Control (CDC)}.\hskip 1em plus 0.5em minus
  0.4em\relax IEEE, 2016, pp. 4252--4272.

\bibitem{fiedler_algebraic_1973}
M.~Fiedler, ``Algebraic connectivity of graphs,'' vol.~23.

\bibitem{olfati04}
R.~{Olfati-Saber} and R.~M. {Murray}, ``Consensus problems in networks of
  agents with switching topology and time-delays,'' \emph{IEEE Transactions on
  Automatic Control}, vol.~49, no.~9, pp. 1520--1533, 2004.

\bibitem{ren07}
W.~Ren and E.~Atkins, ``Distributed multi-vehicle coordinated control via local
  information exchange,'' \emph{International Journal of Robust and Nonlinear
  Control}, vol.~17, pp. 1002--1033, 2007.

\bibitem{gennaro06}
M.~C. De~Gennaro and A.~Jadbabaie, ``Decentralized control of connectivity for
  multi-agent systems,'' in \emph{Proceedings of the 45th IEEE Conference on
  Decision and Control}, 2006, pp. 3628--3633.

\bibitem{nedic18}
A.~Nedi{\'{c}}, A.~Olshevsky, and W.~Shi, \emph{Decentralized Consensus
  Optimization and Resource Allocation}, 2018, pp. 247--287.

\bibitem{holohan2018bounded}
N.~Holohan, S.~Antonatos, S.~Braghin, and P.~Mac~Aonghusa, ``The bounded
  laplace mechanism in differential privacy,'' \emph{arXiv preprint
  arXiv:1808.10410}, 2018.

\bibitem{wang13}
Y.~Wang, X.~Wu, and L.~Wu, ``Differential privacy preserving spectral graph
  analysis,'' in \emph{Pacific-Asia Conference on Knowledge Discovery and Data
  Mining}, 2013, pp. 329--340.

\bibitem{chen2021edge}
B.~Chen, C.~Hawkins, K.~Yazdani, and M.~Hale, ``Edge differential privacy for
  algebraic connectivity of graphs,'' in \emph{2021 60th IEEE Conference on
  Decision and Control (CDC)}.\hskip 1em plus 0.5em minus 0.4em\relax IEEE,
  2021, pp. 2764--2769.

\bibitem{fiedler1975property}
M.~Fiedler, ``A property of eigenvectors of nonnegative symmetric matrices and
  its application to graph theory,'' \emph{Czechoslovak Mathematical Journal},
  vol.~25, no.~4, pp. 619--633, 1975.

\bibitem{Mesbahi2010}
M.~Mesbahi and M.~Egerstedt, \emph{Graph Theoretic Methods in Multiagent
  Networks}, 2010.

\bibitem{gohari21}
P.~{Gohari}, B.~{Wu}, C.~{Hawkins}, M.~{Hale}, and U.~{Topcu}, ``Differential
  privacy on the unit simplex via the dirichlet mechanism,'' \emph{IEEE
  Transactions on Information Forensics and Security}, vol.~16, pp. 2326--2340,
  2021.

\bibitem{PALDINO2017201}
M.~J. Paldino, W.~Zhang, Z.~D. Chu, and F.~Golriz, ``Metrics of brain network
  architecture capture the impact of disease in children with epilepsy,''
  \emph{NeuroImage: Clinical}, vol.~13, pp. 201--208, 2017.

\bibitem{Mohar1991eigenvalues}
B.~Mohar, ``Eigenvalues, diameter, and mean distance in graphs,'' \emph{Graph.
  Comb.}, 1991.

\bibitem{olfati2007consensus}
R.~Olfati-Saber, J.~A. Fax, and R.~M. Murray, ``Consensus and cooperation in
  networked multi-agent systems,'' \emph{Proceedings of the IEEE}, vol.~95,
  no.~1, pp. 215--233, 2007.

\bibitem{jadbabaie2018scaling}
A.~Jadbabaie and A.~Olshevsky, ``Scaling laws for consensus protocols subject
  to noise,'' \emph{IEEE Transactions on Automatic Control}, vol.~64, no.~4,
  pp. 1389--1402, 2018.

\bibitem{levene2002kemeny}
M.~Levene and G.~Loizou, ``Kemeny's constant and the random surfer,'' \emph{The
  American mathematical monthly}, vol. 109, no.~8, pp. 741--745, 2002.

\bibitem{mohar1989isoperimetric}
B.~Mohar, ``Isoperimetric numbers of graphs,'' \emph{Journal of combinatorial
  theory, Series B}, vol.~47, no.~3, pp. 274--291, 1989.

\bibitem{godsil2001algebraic}
C.~Godsil and G.~F. Royle, \emph{Algebraic graph theory}.\hskip 1em plus 0.5em
  minus 0.4em\relax Springer Science \& Business Media, 2001, vol. 207.

\bibitem{bernstein2009matrix}
D.~S. Bernstein, \emph{Matrix mathematics: theory, facts, and formulas}.\hskip
  1em plus 0.5em minus 0.4em\relax Princeton university press, 2009.

\bibitem{kirkland2010algebraic}
S.~Kirkland, ``Algebraic connectivity for vertex-deleted subgraphs, and a
  notion of vertex centrality,'' \emph{Discrete Mathematics}, vol. 310, no.~4,
  pp. 911--921, 2010.

\end{thebibliography}
